\documentclass[12pt, reqno]{amsart}
\usepackage{pifont}
\usepackage{mathrsfs}
\usepackage{geometry}
\usepackage{titletoc}
\usepackage{stix2}
\usepackage{amscd}

\usepackage{amsmath}
\usepackage{amssymb} 
\usepackage{enumitem} 
\usepackage{mathtools} 
\usepackage[table]{xcolor} 
\usepackage[all]{xy} 
\usepackage{tikz} 
\usepackage{tikz-cd}
\usepackage{indentfirst} 
\usepackage{babel} 
\usepackage{setspace} 

\usepackage[colorlinks,linkcolor=red,anchorcolor=green,citecolor=blue]{hyperref} 
\hypersetup{linktocpage = true} 

\usepackage{rotating} 

\usepackage{ytableau} 
\usepackage{longtable} 
\newcolumntype{M}[1]{>{\centering\arraybackslash}m{#1}} 

\usepackage{fancyhdr}

\newcommand\sE{{\mathscr E}}

\newcommand\bp{{\bar\partial}}

\usepackage{amsthm}
\theoremstyle{plain}

\newtheorem{thm}{Theorem}[section]
\newtheorem{lemma}[thm]{Lemma}
\newtheorem{prop}[thm]{Proposition}
\newtheorem{cor}[thm]{Corollary}
\newtheorem{defn}[thm]{Definition}

\theoremstyle{definition}
\newtheorem{example}[thm]{Example}
\newtheorem{remark}[thm]{Remark}

\newcommand{\btheorem}{\begin{thm}}
	\newcommand{\etheorem}{\end{thm}}
\newcommand{\bproposition}{\begin{prop}}
	\newcommand{\eproposition}{\end{prop}}
\newcommand{\bdefinition}{\begin{defn}}
	\newcommand{\edefinition}{\end{defn}}
\newcommand{\bcorollary}{\begin{cor}}
	\newcommand{\ecorollary}{\end{cor}}
\newcommand{\bproof}{\begin{proof}}
	\newcommand{\eproof}{\end{proof}}
\newcommand{\bremark}{\begin{remark}}
	\newcommand{\eremark}{\end{remark}}
\newcommand{\eexample}{\end{example}}
\newcommand{\bexample}{\begin{example}}

\newcommand{\elemma}{\end{lemma}}
\newcommand{\blemma}{\begin{lemma}}

\newcommand{\la}{\langle}
\newcommand{\ra}{\rangle}
\newcommand{\sq}{\sqrt{-1}}

\newcommand{\p}{\partial}

\renewcommand{\bar}{\overline}
\newcommand{\eps}{\varepsilon}

\renewcommand{\phi}{\varphi}

\newcommand{\beq}{\begin{equation}}
\newcommand{\eeq}{\end{equation}}
\newcommand{\ee}{\end{eqnarray*}}
\newcommand{\be}{\begin{eqnarray*}}

\newcommand{\bd}{\begin{enumerate}}
	\newcommand{\ed}{\end{enumerate}}

\renewcommand{\tilde}{\widetilde}

\newcommand{\qtq}[1]{\quad\mbox{#1}\quad}
\renewcommand{\bp}{\bar{\partial}}
\newcommand{\Om}{\Omega}

\newcommand{\ts}{\otimes}

\renewcommand{\>}{\rightarrow}

\newcommand{\C}{{\mathbb C}}

\newcommand{\R}{{\mathbb R}}

\renewcommand{\>}{\rightarrow}

\renewcommand{\p}{{\partial}}
\renewcommand{\bp}{{\bar{\partial}}}

\renewcommand{\bar}{\overline}
\renewcommand{\tilde}{\widetilde}

\newcommand{\smo}{\sqrt{-1}}

\newcommand{\nm}[1]{\left\Vert #1\right\Vert}
\setlist[itemize]{leftmargin=*}
\setlist[enumerate]{leftmargin=*}

\numberwithin{equation}{section} 

\makeatletter

\usepackage{fancyhdr}
\title{Existence of Hermitian metrics  with  prescribed Hermitian-Yang-Mills tensors  I}

\author{Mingwei Wang}
\author{Xiaokui Yang}
\author{Shing-Tung Yau}

\address{Mingwei Wang, Qiuzhen College, Tsinghua University, Beijing, 100084, China}
\email{wangmw21@mails.tsinghua.edu.cn}

\address{Xiaokui Yang, Department of Mathematics and Yau Mathematical Sciences Center, Tsinghua University, Beijing, 100084, China}
\email{xkyang@mail.tsinghua.edu.cn}

\address{Shing-Tung Yau, Yau Mathematical Sciences Center and  Qiuzhen College, Tsinghua University, Beijing, 100084, China}
\email{styau@mail.tsinghua.edu.cn}

\begin{document}
	
	\begin{abstract} In this paper, we solve the prescribed Hermitian-Yang-Mills tensor problem. 
		Let $ E $ be a holomorphic vector bundle over a compact K\"ahler manifold $(M,\omega_g) $. Suppose that there exists a Hermitian metric $ h_0 $ on $E$ such that  its Hermitian-Yang-Mills tensor $ \Lambda_{\omega_g}\left(\sq  R^{h_0}\right) \in \Gamma(M,E^*\otimes \bar E^*) $ is  positive-definite. Then for any positive-definite  Hermitian tensor $ P\in \Gamma(M,E^*\otimes \bar E^*) $, there exists  a unique  Hermitian metric $ h $ on $E$ such that  $$\Lambda_{\omega_g} \left(\sq R^h\right)=P.$$
 The proof is based on a new comparison principle for Hermitian-Yang-Mills tensors. Inspired by these results, we have also derived quantitative Chern number inequalities that apply to both holomorphic vector bundles and compact K\"ahler manifolds.
		
	\end{abstract}
	
	\maketitle {
	
	\setcounter{tocdepth}{1}

	{\small{    \begin{spacing}{1.1} \tableofcontents %
				\dottedcontents{section}[1.8cm]{}{3em}{5pt} %
\end{spacing} }} }
	
	\section{Introduction}
	
	Let $(M,\omega_g)$ be a compact K\"ahler manifold. It is well-known that the class of Ricci curvature $\mathrm{Ric}(\omega_g)$ represents the first Chern class $c_1(M)$ up to a factor of $2\pi$, i.e., 
	$$\left[\mathrm{Ric}(\omega_g)\right]=2\pi c_1(M)\in H^{1,1}(M,\R).$$
	 The converse--given a cohomology class representative, can one find a K\"ahler metric whose Ricci form matches it--is the celebrated content of Yau's solution to the Calabi conjecture \cite{Cal57,Yau78}. Specifically, for any closed real $(1,1)$ form $\Om$ representing $2\pi c_1(M)$,  there exists a unique K\"ahler metric $\omega\in [\omega_g]$ such that 
	\beq \mathrm{Ric}(\omega)=\Om. \label{PRicci}\eeq 
	By using $\p\bp$-lemma, this equation is reduced to a Monge-Amp\`ere type equation 
	$$\left(\omega_g+\sq\p\bp\phi\right)^n=e^F \omega_g^n.$$ 
The Calabi-Yau theorem has indeed been extended beyond K\"ahler manifolds to more general Hermitian manifolds, particularly through the work of G. Sz\'ekelyhidi, V.  Tosatti and B. Weinkove,  among others.  We refer to \cite{TW10a, TW10b, STW17} and the references therein.\\

	The Donaldson-Uhlenbeck-Yau theorem (\cite{UY86}, \cite{Don85, Don87}) asserts  that for  a stable holomorphic vector bundle $E$ over a compact K\"ahler manifold $(M,\omega_g)$,  there exists a unique Hermitian-Einstein metric $h$ on $E$ up to scaling,  satisfying
	\beq \Lambda_{\omega_g}\left(\sq R^h\right)=\lambda_E\cdot h, \label{HE} \eeq 
where $R^h\in\Gamma(M,\Lambda^{1,1}T^*M\ts E^*\ts \bar E^*)$ is the Chern curvature tensor of $(E,h)$ and
 $$\lambda_E=\frac{2\pi n\int_M c_1(E)\wedge \omega_g^{n-1}}{\mathrm{rk}(E)\int_M\omega_g^n}.$$ 
 In particular, the Chern number inequality holds	$$ \int_M \left( (r-1)c^2_1(E)-2rc_2(E) \right) \wedge \omega_g^{n-2} \leq 0. $$
 Together with  results from \cite{Kob82, Lub83}, this establishes a precise correspondence between the algebro-geometric notion of stability and the existence of Hermitian-Einstein metrics. This was generalized to the setting of Hermitian manifolds through the  work \cite{LY86} of  Jun Li and S.-T. Yau.\\
	
The main result of this paper is the following analogue of the Calabi-Yau theorem, addressing the prescribed Hermitian-Yang-Mills tensor problem for holomorphic vector bundles on compact K\"ahler manifolds.
		\btheorem\label{main}
	Let $ E $ be a holomorphic vector bundle over a compact K\"ahler manifold $(M,\omega_g) $. Suppose that there exists a  Hermitian metric $ h_0 $ on $E$ such that its Hermitian-Yang-Mills tensor $ \Lambda_{\omega_g}\left(\sq  R^{h_0}\right) $ is  positive-definite. Then for any  positive-definite Hermitian tensor $ P\in \Gamma(M,E^*\otimes \bar E^*) $,  there exists  a unique  Hermitian metric $ h $ on $E$ such that  \beq  \Lambda_{\omega_g}\left(\sq  R^{h}\right)=P. \label{PHYM}\eeq 
	\etheorem

	\noindent The positivity of  $ \Lambda_{\omega_g}\left(\sq  R^{h_0}\right)$ can be weakened to the positivity of the integral of the minimum eigenvalue of
	 $ \Lambda_{\omega_g}\left(\sq  R^{h_0}\right) $ over the manifold (Theorem \ref{main4}). 
	Moreover, by Serre duality, the corresponding conjugate equation of \eqref{PHYM} admits a solution when $  \Lambda_{\omega_g}\left(\sq  R^{h_0}\right) $ is negative definite; see  Remark \ref{dual} for further details.  The method for solving equation \eqref{PHYM} differs significantly from those for the Calabi-Yau equation \eqref{PRicci} and the Hermitian-Einstein \eqref{HE}, and we provide a brief overview of the key ideas involved. The uniqueness part of Theorem \ref{main} is established by the following comparison principle for Hermitian-Yang-Mills tensors of Hermitian vector bundles:
	
	\btheorem\label{main2}
	Let $E$ be a holomorphic vector bundle over a compact K\"ahler manifold $(M,\omega_g)$.  If $h$ and $h_0$ are  Hermitian metrics on $E$ satisfying  $ \Lambda_{\omega_g}\left(\sq  R^{h_0}\right)>0 $ and \beq   \Lambda_{\omega_g}\left(\sq  R^{h}\right)\leq  \Lambda_{\omega_g}\left(\sq  R^{h_0} \right), \label{Hermitiantensor}\eeq  as Hermitian tensors in $\Gamma(M,E^*\otimes \bar E^*)$, then $ h \leq h_0 $.
	\etheorem
	
Note that inequality \eqref{Hermitiantensor} holds as Hermitian tensors in $\Gamma(M,E^*\otimes  \bar E^*)$,  \textbf{but not  as scaling invariant endomorphisms} in  $\Gamma(M,E^*\otimes E)$. \bremark  The curvature constraint $  \Lambda_{\omega_g}\left(\sq  R^{h_0}\right)>0$ in Theorem \ref{main} is closely related to the notions of RC-positivity and uniform RC-positivity introduced in \cite{Yang18, Yang20}.  A comparison theorem analogous to Theorem \ref{main2} for RC-positive vector bundles was also established in \cite{XYY24+}.  This curvature constraint should be replaced by RC-positivity, and this replacement is carried out in a separate project on general Hermitian manifolds. \eremark

\noindent 	To establish the existence result, we first introduce two key function spaces: let $\mathrm{Herm}(E)$ denote the space of smooth Hermitian tensors on $E$, and  $\mathrm{Herm}^+(E)$ be the space of  Hermitian  metrics on $E$. We then define a \emph{Hermitian-Yang-Mills map}  $$G: \mathrm{Herm}^+(E)\> \mathrm{Herm}(E)$$ by the formula
\beq  G(h)=\Lambda_{\omega_g}\left(\sq  R^{h}\right). \eeq 
We show that if  there exists a  Hermitian metric $ h_0 $ such that $ \Lambda_{\omega_g}\left(\sq  R^{h_0}\right)>0 $, then the restriction map  $$G: G^{-1}\left(\mathrm{Herm}^+(E)\right)\> \mathrm{Herm}^+(E)$$ is a bijection.  The injectivity is given by Theorem \ref{main2}. For the surjectivity,   we show that the image of the restriction map $G$	\beq \mathscr{R} = \{ P \in {\mathrm{Herm}^+(E)} \ |\  G(h)=\Lambda_{\omega_g}\left(\sq  R^{h}\right)=P \text{ for some}\  h \in {\mathrm{Herm}^+(E)} \}\eeq 
	is both open and closed in $\mathrm{Herm}^+(E)$ with respect to the $C^\infty$-topology. The openness of $\mathscr{R}$ (Theorem \ref{SecRicciopen}) is established via the implicit function theorem, leveraging the bijectivity of the linearization of $G$  (up to isomorphisms) at a reference metric $h_1$: \beq \mathscr{L}(\Psi) = \Delta_{\partial^{h_1}} \Psi + \Omega_1 \cdot \Psi,\eeq where the positive definiteness of the Hermitian-Yang-Mills tensor $\Omega_1$ of $h_1$ ensures that $\mathscr{L}$ is an invertible self-adjoint elliptic operator. The closeness of $\mathscr{R}$ (Theorem \ref{Ricciclose}) follows from two key analytic estimates: a scaled version of Theorem \ref{main2}, which provides a uniform $C^0$-bound on metrics $h$ satisfying $G(h) = P$ for $P \in \mathscr{R}$, and a  higher-order regularity estimate (Theorem \ref{ThmC1estimate}), which controls the $C^1$-norm of such metrics. Together, these estimates imply that any sequence $\{P_k\} \subset \mathscr{R}$ converging in the $C^\infty$-topology has a preimage sequence $\{h_k\}$ that is precompact in the $C^1$-norm. By using $W^{k,p}$ and $C^{k,\alpha}$-estimates for elliptic equations, we obtain a smooth limit metric $h$  of some subsequence  $\{h_{k_s}\}$ with $G(h) = \lim P_{k_s}=P$.  Hence $\mathscr{R}$ is closed.\\
	
	      The Calabi-Yau theorem asserts that on a Fano manifold $M$,  there exist K\"ahler metrics $\omega_g $ in each K\"ahler class satisfying $\mathrm{Ric}(\omega_g)>0$. In particular, $\Lambda_{\omega_g} \left(\sq R^g\right)>0$.  As an immediate application of Theorem \ref{main}, we obtain the following result on the existence of Hermitian metrics with prescribed Hermitian-Yang-Mills tensors on Fano manifolds:

	\bcorollary Let $M$ be a Fano manifold. Suppose that  $ P\in \Gamma(M,T^{*1,0}M\otimes  T^{*0,1}M) $ is a  positive-definite  Hermitian section. Then for any K\"ahler class $[\omega]$, there exist a K\"ahler metric $\omega_g\in [\omega]$,  and a unique Hermitian metric $ h $ on $T^{1,0}M$ such that 
	\beq   \Lambda_{\omega_g} \left(\sq R^h\right)=P. \eeq
	\ecorollary
	\noindent 
	The following special case is particularly worth highlighting:
		\bcorollary\label{Fano} On a Fano manifold $M$, for any K\"ahler metric $\omega_g$ with  $\mathrm{Ric}(\omega_g)>0$, there exists a unique Hermitian metric $h $ on $T^{1,0}M$ such that 
	\beq   \Lambda_{\omega_g} \left(\sq R^{ h}\right)=g. \eeq
	Moreover, if $\mathrm{Ric}(\omega_g)=c_0^{-1}\omega_g$ for some $c_0>0$, then $h=c_0g$.
	\ecorollary
	
As inspired by Theorem \ref{main} and the comparsion result described in Theorem \ref{main2}, we derive a quantitative version of Chern number inequalities, extending those associated with Hermitian-Einstein metrics.
	\btheorem\label{main3} Let $(M,\omega_g)$ be a compact K\"ahler manifold and $E$ be a holomorphic vector bundle of rank $r$ over $M$. Suppose that there exists a Hermitian metric $h$ on $E$ satisfying
	\beq  a\cdot h \leq  \Lambda_{\omega_g} \left(\smo  R^{h}\right) \leq b \cdot h, \label{almoststable}\eeq 
	for some  constants $a,b\in \R$.  Then the following Chern number inequality holds
	\beq \int_M \left( (r-1)c^2_1(E)-2rc_2(E) \right) \wedge \omega_g^{n-2} \leq \frac{r(r-1)\left(b - a\right)^2}{8\pi^2n^2}\int_M \omega_g^n.\label{CNI}\eeq
	\etheorem
	
	\noindent Moreover, we establish the following Chern number inequalities extending those associated with K\"ahler-Einstein metrics (\cite{Yau77}).
	\btheorem\label{main5}
	Let $ (M,\omega_g) $ be a compact K\"ahler manifold with  dimension $n\geq 2$.  If there exist two constants $ a,b \in \mathbb{R} $ such that
	\beq a\omega_g \leq \mathrm{Ric}(\omega_g) \leq b\omega_g, \eeq
	then the following Chern number inquality holds 
	\beq \int_M \left(nc^2_1(M)-2(n+1)c_2(M) \right) \wedge \omega_g^{n-2} \leq \frac{(n^2-2)(b - a)^2}{8\pi^2n^2}\int_M \omega_g^n. \eeq
	\etheorem

	\bremark  It is evident that Theorem \ref{main} bears profound theoretical connections to the Calabi-Yau theory, the Donaldson-Uhlenbeck-Yau theory, and the Kazdan-Warner theory. A number of related research initiatives are currently  underway:
	
	\bd  \item[$\bullet$] A Higgs bundle-theoretic formulation of the main theorems will be presented in a separate work.
	\item[$\bullet$]  A flow-based framework, constructed via variations of energy functionals, will be elaborated in subsequent sections.
	\item[$\bullet$] A refined notion of $\varepsilon$-stability, closely tied to the criterion in \eqref{almoststable}, is under development.
	\item[$\bullet$] Theorem \ref{main} and Theorem \ref{main2} together provide an alternative proof of the Donaldson-Uhlenbeck-Yau theorem.
	\ed 

	\eremark

	\noindent\textbf{Acknowledgements}. The second named author would like to thank  Bing-Long Chen, Jixiang Fu,  Kefeng Liu, Linlin Sun, Song Sun,  Valentino Tosatti, Weiping Zhang and Xiping Zhu   for inspiring  discussions.

\vskip 2\baselineskip

\section{Background materials}
Let $ (M,\omega_g) $ be a compact K\"ahler manifold of  dimension $n$ and $ E $ be a holomorphic vector bundle over $M$ of rank $r$.  Fix a smooth Hermitian metric $ h_0 $ on $ E $, the Chern connection of $ (E,h_0) $ is denoted by $ \nabla^{h_0} $.  We shall use the natural decomposition
$$ \nabla^{h_0}=\p^{h_0}+\bp, $$
where $\p^{h_0}$ is the $(1,0)$ part and $\bp$ is the $(0,1)$ part.
The formal adjoints of $ \p^{h_0} $ and $ \bp $ with respect to  $ \omega_g$ and $ h_0 $ are denoted by $ \p^{*,h_0} $ and $ \bar\p{}^{*,h_0} $ respectively. Let $ R^{h_0}$  be the Chern curvature tensor of $ (E,h_0) $.  In local holomorphic coordinates $\{z^i\}$ of $M$ and local holomorphic basis $ \{ e_\alpha \} $ of $E$,  the Chern curvature tensor $R^{h_0}$ can be written as
\beq R^{h_0}= R^{h_0}_{i\bar j\alpha\bar \beta} dz^i\wedge d\bar z^j \ts e^\alpha\ts e_\beta\in\Gamma(M,\Lambda^{1,1}T^*M\ts E^*\ts \bar E^*).\eeq 
The \textbf{Hermitian-Yang-Mills tensor} $ S^{h_0}\in
\Gamma(M,E^*\ts \bar E^*)$ of $(E,h_0)$  is defined as \beq S^{h_0}
:=\Lambda_{\omega_g}\left(\smo R^{h_0}\right)=\left(g^{i\bar
	j}R^{h_0}_{i\bar j\alpha\bar\beta}\right) e^\alpha\ts \bar e^\beta \in
\Gamma(M,E^*\ts \bar E^*).
\eeq  
We say that $S^{h_0}$ is positive-definite, $S^{h_0}>0$, if $\left(S^{h_0}_{\alpha\bar\beta}\right)=\left(g^{i\bar
	j}R^{h_0}_{i\bar j\alpha\bar\beta}\right)$ is a positive-definite Hermitian matrix at each point of $M$.
By using the Hermitian metric $h_0$, there
are natural lifts of $R^{h_0}$ and $S^{h_0}$: \beq\Theta^{h_0}=R^{h_0}\cdot
h_0^{-1}=\left(R^{h_0}\right)_{i\bar j\alpha}^{\beta}dz^i\wedge d\bar z^j\ts
e^\alpha\ts{e}_\beta\in \Gamma(M,\Lambda^{1,1}T^*M\ts E^*\ts E),\eeq
and \beq   K^{h_0}=S^{h_0}\cdot h_0^{-1}=\left(g^{i\bar j}\left(R^{h_0}\right)_{i\bar j\alpha}^{\beta}\right) e^\alpha\ts  e_\beta\in \Gamma(M,E^*\ts E).
\eeq
The Chern scalar curvature $s_{h_0}$ is defined as 
\beq s_{h_0}=\mathrm{tr}_E K^{h_0}=\mathrm{tr}_{h_0} S^{h_0}=g^{i\bar j}h_0^{\alpha\bar\beta}R^{h_0}_{i\bar j \alpha\bar\beta}.\eeq 

 Let $h\in\Gamma(M,E^*\ts{\bar E}^*)$ be another smooth Hermitian metric on $E$. 
  There is an induced smooth section $ H \in \Gamma(M,E^* \otimes E) $ given by 
 \beq H = h \cdot h_0^{-1}.\eeq 
In local holomorphic basis $ \{ e_\alpha \} $ of $E$,  $H=H_\alpha^\beta e^\alpha\ts e_\beta$ is  represented by 
\beq H_\alpha^\beta = h_{\alpha\bar\gamma} h_0^{\beta\bar\gamma}. \eeq
It is clear that the inverse $ H^{-1}=(H^{-1})_\alpha^\beta e^\alpha\ts e_\beta \in \Gamma(M,E^*\otimes E) $ has components 
\beq (H^{-1})_\alpha^\beta =h_{0,\alpha\bar\gamma} h^{\beta\bar\gamma}. \eeq
For $ P=P_\alpha^\beta e^\alpha\ts e_\beta \in \Omega^p(M,E^*\otimes E) $ and $ Q=Q_\alpha^\beta e^\alpha\ts e_\beta \in \Omega^q(M,E^*\otimes E) $, the product \beq  P\cdot Q \in \Omega^{p+q}(M,E^*\otimes E) \eeq  is defined as
\beq P \cdot Q = (P_\alpha^\beta \wedge Q_\beta^\gamma) \; e^\alpha \otimes e_\gamma. \eeq 
In particular, one has \beq  H \cdot H^{-1} = H^{-1} \cdot H = \mathrm{Id}_E \in \Gamma(M,E^*\otimes E). \eeq

\noindent The tensor $H$ will play a central role in what follows; we now present some of its key linear algebraic properties.

\bdefinition  A  section $P \in\Gamma(M,E^*\ts E)$ is called $h_0$-Hermitian if 
\beq h_0(Pv,w)=h_0(v,Pw), \eeq 
for any $ v, w\in \Gamma(M,E)$.  It is easy to see that $P=P_\alpha^\beta e^\alpha\ts e_\beta$ is $h_0$-Hermitian if and only if 
\beq h_{0,\alpha\bar\delta} h_0^{\beta\bar\gamma} \bar P{}_\gamma^\delta =P_\alpha^\beta.\eeq 
An $h_0$-Hermitian section $P=P_\alpha^\beta e^\alpha\ts e_\beta \in\Gamma(M,E^*\ts E)$ is  called positive-definite if $$h_0(Pv,v)>0$$ for all $v\neq 0$. This is equivalent to the  Hermitian  matrix
$\left(P_{\alpha\bar \beta}\right)=\left(h_{0,\gamma\bar\beta}P_{\alpha}^\gamma\right)$ being positive-definite. We also write  $P>0$. Similarly, $P\geq 0$ is defined by requiring   $h_0(Pv,v)\geq 0$.
\edefinition

\blemma \label{linearalgebra}\bd 
 \item $H=h\cdot h_0^{-1}\in\Gamma(M,E^*\ts E)$    is  $h_0$-Hermitian.

\item If $P\in \Gamma(M,E^*\ts E)$ is $h_0$-Hermitian,  for any $ A, B \in \Gamma(M,E^*\otimes E) $,
\beq h_0(P\cdot A, B)=h_0(A, P\cdot B), \quad  h_0(A \cdot P, B) = h_0(A, B\cdot P). \eeq 

\item If $ A, B \in \Gamma(M,E^*\otimes E) $ are $h_0$-Hermitian and $ A\geq 0, B \geq 0 $, then  
\beq h_0(A,B) \geq 0. \eeq Moreover,  for  any section $ C \in \Gamma(M,E^* \otimes E) $, one has
\beq h_0(A \cdot C \cdot B, C ) \geq 0. \eeq

\ed
\elemma

\bproof 
 (1) follows from a straightforward computation that
\beq h_{0,\alpha\bar\delta} h_0^{\beta\bar\gamma} \bar H{}_\gamma^\delta = h_{0,\alpha\bar\delta} h_0^{\beta\bar\gamma} h_0^{p\bar\delta} h_{p\bar\gamma} =  H_\alpha^\beta. \eeq
For (2), one has
\be h_0(A\cdot P, B) & = & h_{0,\gamma\bar q}h_0^{\alpha\bar p}A_\alpha^\beta P_\beta^\gamma \bar B{}_p^q =  h_{0,\gamma\bar q}h_0^{\alpha\bar p}A_\alpha^\beta \left( h_{0,\beta\bar s} h_0^{\gamma\bar t} \bar P{}_t^s\right) \bar B{}_p^q \\
& = & h_{0,\beta\bar s}h_0^{\alpha\bar p}A_\alpha^\beta \bar B{}_p^q \bar P{}_q^s = h_0(A,B\cdot P), \ee
where the second identity holds since $ P $ is $h_0$-Hermitian. The relation $h_0(P\cdot A, B)=h_0(A, P\cdot B)$ can be verified in a similar way.\\

\noindent  For (3), since $A\geq 0$ and $B\geq 0$, $A^{1/2}$ and $B^{1/2}$ are well-defined and it is easy to deduce that both of them are $h_0$-Hermitian. Hence,
\beq h_0(A,B)=h_0\left(A^{1/2},  A^{1/2}B\right)=h_0\left(A^{1/2}B^{1/2},  A^{1/2}B^{1/2}\right)\geq 0. \eeq 

\noindent A simple  calculation shows that
$h_0(A \cdot C \cdot B, C )= h_0(A^{1/2} \cdot C \cdot B^{1/2}, A^{1/2} \cdot C \cdot B^{1/2} ) \geq 0$.
\eproof

\vskip 2\baselineskip

\section{A comparison theorem for Hermitian-Yang-Mills tensors}

In this section, we present a proof of Theorem \ref{main2}, which provides a comparison result for Hermitian-Yang-Mills tensors of abstract vector bundles. We then explore its applications, including the uniqueness of solutions to the prescribed Hermitian-Yang-Mills tensor equation and a version of
$C^0$-estimate for Hermitian-Yang-Mills tensors.
\btheorem	\label{CompareThm} 
	Let $E$ be a holomorphic vector bundle over a compact K\"ahler manifold $(M,\omega_g)$.  If $h$ and $h_0$ are  Hermitian metrics on $E$ satisfying  $ \Lambda_{\omega_g}\left(\sq  R^{h_0}\right)>0 $ and \beq   \Lambda_{\omega_g}\left(\sq  R^{h}\right)\leq  \Lambda_{\omega_g}\left(\sq  R^{h_0} \right), \eeq  as Hermitian tensors in $\Gamma(M,E^*\otimes \bar E^*)$, then $ h \leq h_0 $.
\etheorem

\noindent 
As an application of Theorem \ref{CompareThm}, one has 

\bcorollary	\label{uniqueness} 
Let $E$ be a holomorphic vector bundle over a compact K\"ahler manifold $(M,\omega_g)$. Suppose that $h$ and $h_0$ are  Hermitian metrics on $E$.	If \beq   \Lambda_{\omega_g}\left(\sq   R^{h}\right)=  \Lambda_{\omega_g}\left(\sq  R^{h_0}\right)>0 \eeq  as Hermitian tensors in $\Gamma(M,E^*\otimes \bar E^*)$, then $ h = h_0 $.
\ecorollary

\noindent In particular, we establish the uniqueness in Theorem \ref{main}:

\bcorollary
Let $ E $ be a holomorphic vector bundle over a compact K\"ahler manifold $(M,\omega_g) $. Then for any Hermitian positive definite section $ P \in \Gamma(M,E^*\otimes \bar E^*) $, there exists  at most one  Hermitian metric $ h $ on $E$ such that
\beq \Lambda_{\omega_g} \left(\sq R^h\right)=P.\eeq 
\ecorollary

\noindent By rescaling, we obtain the following result, which is analogous to the $C^2$-estimate in the proof of the Calabi-Yau Theorem.
\bcorollary\label{scaledcompare} Suppose that $h$ and $h_0$ are two smooth Hermitian metrics on $E$. If  $ \Lambda_{\omega_g}\left(\sq  R^{h_0}\right)>0 $ and \beq    \Lambda_{\omega_g}\left(\sq  R^{h}\right)\leq \lambda\cdot   \Lambda_{\omega_g}\left(\sq  R^{h_0}\right) \eeq  as Hermitian tensors in $\Gamma(M,E^*\otimes \bar E^*)$ and $\lambda\in \R^+$, then $ h \leq \lambda \cdot h_0 $.

\ecorollary

\noindent The proof of Theorem \ref{CompareThm} starts with the following relation concerning  Hermitian-Yang-Mills tensors.
\bproposition 	Let $E$ be a holomorphic vector bundle over a compact K\"ahler manifold $(M,\omega_g)$. Suppose that $h$ and $h_0$ are two  Hermitian metrics on $E$. Then \beq K^{h}-K^{h_0}= \Lambda_{\omega_g}\left(\sq  \Theta^{h}\right)- \Lambda_{\omega_g}\left(\sq  \Theta^{h_0}\right)= \p^{*,h_0}\left( (\p^{h_0} H) \cdot H^{-1} \right), \label{secondChernRiccidiff} \eeq
as tensors in $\Gamma(M,E^*\ts E)$, where  $H=h\cdot h_0^{-1}$ and $\p^{h_0}$ is  the $(1,0)$ component of the Chern connection on $E^*\ts E$ induced by $(E,h_0)$. 
\eproposition

\bproof The Chern connection $\nabla^{h_0}$ is given by
\beq \nabla_i^{h_0}e_\beta =(\Gamma_0)_{i\beta}^\mu e_\mu, \quad \nabla_i^{h_0}e^\alpha = - (\Gamma_0)_{i\beta}^\alpha e^\beta.  \eeq
where $(\Gamma_0)_{i\alpha}^\beta$ is the Christoffel symbol of the Chern connection on $(E,h_0)$. Hence, 
\be
\nabla_i^{h_0} H 
&=&  \nabla^{h_0}_i\left(h_0^{\beta\bar\gamma} h_{\alpha\bar\gamma} \; e^\alpha \otimes e_\beta \right) \\
& = & \p_i\left(  h_0^{\beta\bar\gamma} h_{\alpha\bar\gamma}\right) \; e^\alpha \otimes e_\beta +h_0^{\beta\bar\gamma} h_{\alpha\bar\gamma} \left(\nabla_i^{h_0} e^\alpha\right) \otimes e_\beta + h_0^{\beta\bar\gamma} h_{\alpha\bar\gamma} e^\alpha\ts \left(\nabla_i^{h_0} e_\beta\right)   \\
& = & \left(\left(\p_i h_0^{\beta\bar\gamma}\right) h_{\alpha\bar\gamma} +  h_0^{\beta\bar\gamma} (\p_ih_{\alpha\bar\gamma}) - h_0^{\beta\bar\gamma} h_{\mu\bar\gamma} (\Gamma_0)_{i\alpha}^\mu+ h_0^{\mu\bar\gamma} h_{\alpha\bar\gamma} (\Gamma_0)_{i\mu}^\beta\right) \; e^\alpha \otimes e_\beta \\
& = & h_0^{\beta\bar\gamma} \left((\p_ih_{\alpha\bar\gamma}) - h_{\mu\bar\gamma} (\Gamma_0)_{i\alpha}^\mu\right) \; e^\alpha \otimes e_\beta.
\ee
Therefore,  
\beq \p^{h_0} H = h_0^{\beta\bar\gamma} \left((\p_ih_{\alpha\bar\gamma}) -  h_{\mu\bar\gamma} (\Gamma_0)_{i\alpha}^\mu\right) \; (dz^i \otimes e^\alpha \otimes e_\beta) \eeq
and 
\begin{eqnarray} (\p^{h_0} H ) \cdot H^{-1} 
\nonumber& = & h^{\beta\bar\gamma} \left((\p_ih_{\alpha\bar\gamma}) - h_{\mu\bar\gamma} (\Gamma_0)_{i\alpha}^\mu\right) \; (dz^i \otimes e^\alpha \otimes e_\beta) \\
& = & \left(\Gamma_{i\alpha}^\beta - (\Gamma_0)_{i\alpha}^\beta\right) \; (dz^i \otimes e^\alpha \otimes e_\beta)\label{torsion} \end{eqnarray}
where $\Gamma_{i\alpha}^\beta$ is the Christoffel symbol of the Chern connection $\nabla^h$ of $h$.
Moreover, the Chern curvature components are given by
\beq (R^h)_{i\bar j\alpha}{}^\beta =-\frac{\p \Gamma_{i\alpha}^\beta}{\p\bar z^j}, \quad (R^{h_0})_{i\bar j\alpha}{}^\beta =-\frac{\p (\Gamma_0)_{i\alpha}^\beta}{\p\bar z^j}, \eeq
and this implies the curvature identity
\beq \bar\p{}\left( (\p^{h_0} H) \cdot H^{-1} \right) = \left((R^h)_{i\bar j\alpha}{}^{\beta} - (R^{h_0})_{i\bar j\alpha}{}^{\beta}\right) (dz^i \wedge d\bar z{}^j \otimes e^\alpha \otimes e_\beta) = \Theta^h - \Theta^{h_0}. \eeq
On the other hand, by using the Bochner-Kodaira formula \beq  \p^{*,h_0} = \smo[\Lambda_{\omega_g}, \bar\p],  \eeq one obtains the difference  \eqref{secondChernRiccidiff} of the  Hermitian-Yang-Mills tensors.
\eproof

\vskip 1\baselineskip
	
	\bproof[Proof of Theorem \ref{CompareThm}]
 To simplify notations, we set $\Omega = K^{h_0} \in \Gamma(M,E^*\otimes E)$ and define
	\beq  \Phi := h_0^{\alpha\bar\gamma} S^h_{\beta\bar\gamma} \; e^\beta \otimes e_\alpha \in \Gamma(M,E^*\otimes E). \eeq
	It's easy to see that
	\beq  \Phi = S^h\cdot h_0^{-1}=K^h \cdot H. \eeq On the other hand, by \eqref{secondChernRiccidiff},
	\beq K^h = K^{h_0} + \p^{*,h_0}(\p^{h_0}H \cdot H^{-1}), \eeq
	and in particular, we obtain the identity
	\beq \label{Currelation}\Phi\cdot H^{-1} = \Omega + \p^{*,h_0}(\p^{h_0}H \cdot H^{-1}). \eeq
This equation can be rewritten in the following form:
	\beq \label{mainidentity2} \Omega \cdot (H^{-1} - \mathrm{Id}_E) = \p^{*,h_0}(\p^{h_0}H \cdot H^{-1}) + (\Omega - \Phi) \cdot H^{-1}. \eeq
	We claim that \beq H \leq \operatorname{Id}_E \label{keycompare}\eeq  with respect to $h_0$ which also implies $h \leq h_0$. Indeed,  for any $ v = v^\alpha \cdot e_\alpha \in \Gamma(M,E)$,  
	\beq h_0(Hv,v) = h_{0,\beta\bar\gamma} H^\beta_\alpha v^\alpha \bar v{}^\gamma = h_{0,\beta\bar\gamma} h_0^{\beta\bar\delta} h_{\alpha\bar\delta} v^\alpha \bar v{}^\gamma = h_{\alpha\bar\gamma} v^\alpha \bar v{}^\gamma = h(v,v). \eeq
	Hence, the relation $H \leq \operatorname{Id}_E$ implies that for all $v\in \Gamma(M,E)$
	\beq h(v,v)=h_0(Hv,v) \leq h_0(v,v).\eeq 
	
	\noindent
	To prove \eqref{keycompare},  we define $\kappa: M \rightarrow \mathbb{R}$ as
	\beq 
	\kappa(x):=\sup _{v \in E_x, v \neq 0} \frac{ h_0(Hv,v) }{h_0(v,v)} .
\eeq 
The maximal eigenvalue of $H$ with respect to $h_0$ is
	\beq 
	\Lambda:=\sup _{x \in M} \kappa(x) .
\eeq 
We argue by contradiction and assume that \eqref{keycompare} does not hold. In this case, one has \beq  \Lambda > 1. \eeq For any real number $ \eps > 0 $, we consider the following test section
\beq W_\eps:=\left ((\Lambda+\eps)\mathrm{Id}_E - H\right)^{-2n} \cdot H \in \Gamma(M,E^*\ts E).\eeq 
Since the maximal eigenvalue of $H$ is $\Lambda$, 
\beq (\Lambda+\eps)\mathrm{Id}_E - H>0,\eeq 
and so  \beq  \left((\Lambda+\eps)\mathrm{Id}_E - H\right)^{-2n} >0 \eeq  with respect to $h_0$. By part $(2)$ of Lemma \ref{linearalgebra}, one has 
\be \left((\Omega - \Phi) \cdot H^{-1}, W_\eps \right)_{h_0}&=&\left((\Omega - \Phi) \cdot H^{-1}, ((\Lambda+\eps)\mathrm{Id}_E - H)^{-2n} \cdot H \right)_{h_0} \\&=& \left(\Omega - \Phi, ((\Lambda+\eps)\mathrm{Id}_E - H)^{-2n} \right)_{h_0}. \ee
On the other hand, $\Phi$ is  Hermitian with respect to $ h_0 $:
\beq  h_{0,\alpha\bar\delta} h_0^{\beta\bar\gamma} \bar \Phi{}_\gamma^\delta = h_{0,\alpha\bar\delta} h_0^{\beta\bar\gamma} h_0^{p\bar\delta} S^h_{p\bar\gamma} = h_0^{\beta\bar\gamma} S^h_{\alpha\bar\gamma} = \Phi_\alpha^\beta. \eeq
Hence,  by the curvature assumption $ \left(S^h_{\alpha\bar\beta}\right) \leq \left(S^{h_0}_{\alpha\bar\beta}\right)$,  we have \beq  \Omega - \Phi \geq 0 \eeq with respect to $h_0$. By part $(3)$ of Lemma \ref{linearalgebra},
\beq \left((\Omega - \Phi) \cdot H^{-1}, W_\eps \right)_{h_0}= \left(\Omega - \Phi, ((\Lambda+\eps)\mathrm{Id}_E - H)^{-2n} \right)_{h_0}\geq 0.\eeq

\vskip 1\baselineskip

\noindent  For any $ k \geq 0 $ and any fixed point $ x \in M $, we assume that  $ g_{i\bar j}(x) = \delta_{ij} $ and so  at point $x$
\be \langle \p^{h_0} H \cdot H^{-1}, \p^{h_0} H^k \rangle_{h_0} & = & \sum_{p=1}^k \langle \p^{h_0} H \cdot H^{-1}, H^{p-1} \cdot \p^{h_0} H \cdot H^{k-p} \rangle_{h_0} \\ 
& = & \sum_{i=1}^n \sum_{p=1}^k \langle H^{p-1} \cdot \p_i^{h_0} H \cdot H^{k-p-1}, \p_i^{h_0} H \rangle_{h_0}\\& \geq& 0, \ee
where the last inequality holds by part (3) of Lemma \ref{linearalgebra}.  Since $H$ has maximal eigenvalue $\Lambda$, one has the Taylor expansion \beq \label{TaylorH}\left((\Lambda+\eps)\mathrm{Id}_E - H\right)^{-2n} = \sum_{k=0}^{\infty} {2n+k-1\choose k} (\Lambda+\eps)^{-(k+2n)}H^{k}. \eeq
By using this expansion and integration by parts, one gets
\be&& \left(\p^{*,h_0}(\p^{h_0}H \cdot H^{-1}), W_\eps\right)_{h_0}\\&=& \left(\p^{h_0}H \cdot H^{-1}, \p^{h_0}\left(((\Lambda+\eps)\mathrm{Id}_E - H)^{-2n} \cdot H\right)\right)_{h_0} \\
& = & \sum_{k=0}^\infty {2n+k-1\choose k} (\Lambda+\eps)^{-(k+2n)} \left(\p^{h_0}H \cdot H^{-1}, \p^{h_0}H^{k+1}\right)_{h_0} \geq 0. \ee

\noindent 
By pairing with $ W_\eps $ in \eqref{mainidentity2},  one gets
\beq \left(\Omega \cdot (H^{-1}-\mathrm{Id}_E),W_\eps \right)_{h_0} =\left(\p^{*,h_0}(\p^{h_0}H \cdot H^{-1}), W_\eps\right)_{h_0}+ \left((\Omega - \Phi) \cdot H^{-1}, W_\eps \right)_{h_0}\geq 0. \label{keymaximal}\eeq 

\vskip 1\baselineskip

\noindent In the following, we estimate $\left(\Omega \cdot (H^{-1}-\mathrm{Id}_E),W_\eps \right)_{h_0}$
using a different method and derive a contradiction.
Let $ K $ and $ k $ be the largest and smallest eigenvalues of $ \Omega $ with respect to $h_0$ respectively, i.e., 
	\beq K = \sup_{x\in M} \sup_{v \in E_x, v\neq 0}  \frac{ h_0(\Om v,v) }{h_0(v,v)} > 0, \quad k = \inf_{x\in M} \inf_{v \in E_x, v\neq 0}  \frac{ h_0(\Om v,v) }{h_0(v,v)} > 0. \eeq
	In particular, 
	\beq \mathrm{tr}_E\Om\leq r K.\eeq
	where $r$ is the rank of $E$.	We claim that 
	\bd
	\item[(I)] For any point $ x \in M $, one obtains the uniform estimate
	\beq  \left\langle \Omega,((\Lambda+\eps)\mathrm{Id}_E - H)^{-2n} \cdot (\mathrm{Id}_E-H) \right\rangle_{h_0}(x) \leq rK(\Lambda+\eps-1)^{-2n}. \label{I}\eeq
	\item[(II)] For any $ x \in M $ satisfying $ \kappa(x) \geq (\Lambda+1)/2 $,  the following refined estimate holds:
	\beq \left\langle \Omega,((\Lambda+\eps)\mathrm{Id}_E - H)^{-2n} \cdot (\mathrm{Id}_E-H) \right\rangle_{h_0}(x) \leq rK(\Lambda+\eps-1)^{-2n} - \frac{k(\Lambda-1)}{2(\Lambda+\eps-\kappa(x))^{2n}}. \label{II}\eeq
	\ed 
To prove (I), 	let $ f $ be an auxiliary function defined by  
	\beq f(t ) = \frac{1-t}{(\Lambda+\eps-t)^{2n}},\quad t\in(0,\Lambda]. \eeq 
	One can see clearly that $f(H)$ is a well-defined smooth section in $\Gamma(M,E^*\ts E)$ and it is $h_0$-Hermitian.  Moreover, 
	\beq  \left\langle \Omega,((\Lambda+\eps)\mathrm{Id}_E - H)^{-2n} \cdot (\mathrm{Id}_E-H) \right\rangle_{h_0} = \left\langle \Omega, f(H) \right\rangle_{h_0} = \mathrm{tr}_E(\Omega \cdot f(H)), \eeq
	where the last identity holds since both $\Om$ and $f(H)$ are $h_0$-Hermitian.
	Let $$ 0 < \lambda_1 \leq \cdots \leq \lambda_r \leq \Lambda $$ be  eigenvalues of $ H $ at point $ x \in M $. The corresponding  eigenvalues of $ f(H) $ are $$ f(\lambda_1),\cdots , f(\lambda_r).$$ It is clear that  for any $t\in (0,\Lambda]$
	\beq f(t) \leq (\Lambda+\eps-1)^{-2n}.  \eeq
	In particular, one has
	\beq f(H) \leq (\Lambda+\eps-1)^{-2n}\mathrm{Id}_E. \label{operatorinequality}\eeq
	Therefore,  
	\beq \mathrm{tr}_E(\Omega \cdot f(H)) \leq (\Lambda+\eps-1)^{-2n}\mathrm{tr}_E(\Omega) \leq rK(\Lambda+\eps-1)^{-2n}. \eeq
This establishes \eqref{I}. \\

\noindent  For (II),  if  $ x \in M $ satisfies $ \kappa(x) \geq (\Lambda+1)/2 $, then
	\beq f(\lambda_r) = f(\kappa(x)) = \frac{1-\kappa(x)}{(\Lambda+\eps-\kappa(x))^{2n}} \leq -\frac{\Lambda-1}{2(\Lambda+\eps-\kappa(x))^{2n}}. \label{peak}\eeq
Suppose that $ v_x\in E_x$ is  an eigenvector of $H$ with unit length  corresponding to $ \lambda_r = \kappa(x) $. By \eqref{operatorinequality}, \eqref{peak} and the fact that $\Lambda>1$, one achieves the following estimate at point $x$
	\beq f(H) \leq (\Lambda+\eps-1)^{-2n}\mathrm{Id}_E - \frac{\Lambda-1}{2(\Lambda+\eps-\kappa(x))^{2n}} \; v_x^* \otimes v_x . \eeq
In particular, one has 
	\be \mathrm{tr}_E(\Omega \cdot f(H)) & \leq & (\Lambda+\eps-1)^{-2n}\mathrm{tr}_E(\Omega) - \frac{\Lambda-1}{2(\Lambda+\eps-\kappa(x))^{2n}}\mathrm{tr}_E(\Omega \cdot(v_x^*\otimes v_x)) \\
	& \leq & rK(\Lambda+\eps-1)^{-2n} - \frac{k(\Lambda-1)}{2(\Lambda+\eps-\kappa(x))^{2n}}. \ee
	This establishes the inequality \eqref{II}.\\

Now we fix a point $ x_0 \in M $ with $ \kappa(x_0) = \Lambda $. Since $ \kappa(x) $ is a Lipschitz function with  Lipschitz constant $ L=L(H,M,\omega_g, h_0)> 0 $,  we set $$ \delta = \frac{\Lambda-1}{2L} > 0. $$ For any $ x \in M $ with $ d(x,x_0) < \delta $, one has $|\kappa(x)-\kappa(x_0)|\leq L d(x,x_0)$, and so
	\beq \kappa(x) \geq \kappa(x_0) - Ld(x,x_0)=\Lambda - Ld(x,x_0)\geq \frac{\Lambda+1}{2}. \eeq 
	By \eqref{II}, 
	\be && \int_{B(x_0,\delta)}\la \Omega,((\Lambda+\eps)\mathrm{Id}_E - H)^{-2n} \cdot (\mathrm{Id}_E-H) \ra_{h_0} \cdot \frac{\omega_g^n}{n!}\\
	&\leq & \frac{rK}{(\Lambda+\eps-1)^{2n}}\int_{B(x_0,\delta)}\frac{\omega_g^n}{n!} - \int_{B(x_0,\delta)} \frac{k(\Lambda-1)}{2(\Lambda+\eps-\kappa(x))^{2n}} \cdot \frac{\omega_g^n}{n!}\\
	&\leq & \frac{rK}{(\Lambda+\eps-1)^{2n}}\int_{B(x_0,\delta)}\frac{\omega_g^n}{n!}-\int_{B(x_0,\delta)} \frac{k(\Lambda-1)}{2(\eps + Ld(x,x_0))^{2n}} \cdot \frac{\omega_g^n}{n!}  \ee 
	since $\Lambda-\kappa(x)\leq Ld(x,x_0)$.
Moreover, by  \eqref{I}, one obtains 
	\be &&\left( \Omega,((\Lambda+\eps)\mathrm{Id}_E - H)^{-2n} \cdot (\mathrm{Id}_E-H) \right)_{h_0}\\
	& \leq & \frac{rK\mathrm{Vol}(M,\omega_g)}{(\Lambda+\eps-1)^{2n}} - \int_{B(x_0,\delta)} \frac{k(\Lambda-1)}{2(\eps + Ld(x,x_0))^{2n}} \cdot \frac{\omega_g^n}{n!}.  \ee
	On the other hand,  we established in \eqref{keymaximal} that 
	\beq\left( \Omega,((\Lambda+\eps)\mathrm{Id}_E - H)^{-2n} \cdot (\mathrm{Id}_E-H) \right)_{h_0}=\left(\Omega \cdot (H^{-1}-\mathrm{Id}_E),W_\eps \right)_{h_0}\geq 0. \eeq 
By the monotone convergence theorem, one has 
	\be
	0 & \leq & \liminf_{\eps\>0^+} \left( \Omega,((\Lambda+\eps)\mathrm{Id}_E - H)^{-2n} \cdot (\mathrm{Id}_E-H) \right)_{h_0} \\
	& \leq & \frac{rK\mathrm{Vol}(M,\omega_g)}{(\Lambda-1)^{2n}} -  \frac{k(\Lambda-1)}{2L^{2n}}\int_{B(x_0,\delta)} \frac{1}{d(x,x_0)^{2n}}\cdot \frac{\omega_g^n}{n!}=-\infty . \ee
	That's a contradiction. Therefore, one deduces $ \Lambda \leq 1 $ and so $ h \leq h_0 $. \eproof
	
	\vskip 2\baselineskip

	\section{A priori estimates}
	
	In this section, we establish  priori $C^1$-estimates for $H=h\cdot h_0^{-1}$ where $h$ and $h_0$ are two  Hermitian metrics on $E$. Recall that  
	\bd \item $S^h=\Lambda_{\omega_g}\left(\sq R^{h}\right)$ and $S^{h_0}=\Lambda_{\omega_g}\left(\sq R^{h_0}\right)$ are the Hermitian-Yang-Mills  tensors of $h$ and $h_0$ respectively;

	\item  $S^{h_0}=\Lambda_{\omega_g}\left(\sq R^{h_0}\right)>0$;
	
	\item We use notations in $\Gamma(M,E^*\ts E)$
\beq	H =  h\cdot h_0^{-1},\quad \Om=K^{h_0} , \quad \Phi = K^h\cdot H. \eeq 
	\ed
	\noindent The main result of this section is following uniform $C^1$-estimate for $H$. By using this uniform estimate and an elliptic equation  satisfied by $H$, we derive uniform higher-order $W^{k,p}$ and $C^{k,\alpha}$-estimates.
		\btheorem \label{ThmC1estimate} Let $E$ be a holomorphic vector bundle over  a compact K\"ahler manifold $(M,\omega_g)$. Suppose that $h_0$ and $h$ are two Hermitian metrics on $E$, and $\Lambda_{\omega_g}\left(\sq R^{h_0}\right)>0$. 	
		 	If there exists a constant  $C_1\geq 1 $ such that 
	\beq C_1^{-1}h_0 \leq \Lambda_{\omega_g}\left(\sq R^h\right)\leq C_1 h_0 \eeq
	as Hermitian tensors in $\Gamma(M,E^*\ts \bar E^*)$,
	then the following  $C^1$-estimate holds for $H=h\cdot h_0^{-1}$:
	\beq  |H|_{C^1(M,\omega_g, h_0)}\leq C_{16},\eeq
	where $C_{16}$ depends on $M,\omega_g, C_1, h_0$ and an upper bound of  $|\Phi|_{C^1(M,\omega_g, h_0)}$. 
	
	\etheorem
	
\noindent	We begin with the following estimate.
	
	\bproposition
	\label{C0estimate}
Let $E$ be a holomorphic vector bundle over  a compact K\"ahler manifold $(M,\omega_g)$. Suppose that $h_0$ and $h$ are two Hermitian metrics on $E$, and $\Lambda_{\omega_g}\left(\sq R^{h_0}\right)>0$. 	
If there exists a constant  $C_1\geq 1 $ such that 
\beq C_1^{-1}h_0 \leq \Lambda_{\omega_g}\left(\sq R^h\right)\leq C_1 h_0 \label{boundedcurvature}\eeq
as Hermitian tensors in $\Gamma(M,E^*\ts \bar E^*)$,
	then there exists  $ C_2 = C_2(M,\omega_g,h_0,C_1) $ such that 
 \beq C_2^{-1} h_0 \leq h \leq C_2h_0.\eeq 
	\eproposition
	
	\bproof Since $M$ is compact and $\Lambda_{\omega_g}\left(\sq R^{h_0}\right)>0$, there exists $ C_3 = C_3(M,\omega_g,h_0) > 0 $ such that 
	\beq C_3^{-1}h_0 \leq \Lambda_{\omega_g}\left(\sq R^{h_0}\right)\leq C_3 h_0. \eeq
By the condition \eqref{boundedcurvature}, one has
	\beq (C_1C_3)^{-1}\Lambda_{\omega_g}\left(\sq R^{h_0}\right)\leq \Lambda_{\omega_g}\left(\sq R^{h}\right) \leq C_1C_3 \Lambda_{\omega_g}\left(\sq R^{h_0}\right),\eeq
	as Hermitian tensors in $\Gamma(M,E^*\ts \bar E^*)$. By Corollary \ref{scaledcompare}, one has $ C_2^{-1}h_0 \leq h \leq C_2h_0 $ for $ C_2 = C_1C_3$.
	\eproof

	\bproposition \label{W12estimate} Let $E$ be a holomorphic vector bundle over  a compact K\"ahler manifold $(M,\omega_g)$. Suppose that $h_0$ and $h$ are two Hermitian metrics on $E$, and $H=h\cdot h_0^{-1}$.
	 Suppose that there exists a constant $C_4\geq 1$ such that
	\beq C_4^{-1}h_0\leq h\leq C_4h_0. \eeq
	Then for any $A\in\Gamma(M,E^*\ts E)$, 
	\beq C^{-1}_4 h_0(A\cdot H^{-1}, A)\leq h(A\cdot H^{-1}, A\cdot H^{-1}) \leq C_4 h_0(A\cdot H^{-1}, A)\label{estimate1}\eeq 
	and 
	\beq C_4^{-2} h_0(A,A) \leq h(A\cdot H^{-1}, A\cdot H^{-1}) \leq C_4^2 h_0(A,A).\label{c15}\eeq  
In particular,
\beq C_4^{-2} h_0(\p^{h_0} H, \p^{h_0} H) \leq h(\p^{h_0} H\cdot H^{-1}, \p^{h_0} H\cdot H^{-1}) \leq C_4^2 h_0(\p^{h_0} H, \p^{h_0} H). \label{estimate3}\eeq  
Moreover, 
	\beq \mathrm{tr}_E\left(\Lambda_{\omega_g}\left(\sq \p^{h_0} H\cdot H^{-1}\cdot \bp H\right)\right)\geq C_4^{-1} h\left(\p^{h_0}H\cdot H^{-1}, \p^{h_0}H\cdot H^{-1}\right). \label{estimate2}\eeq 
	\eproposition
	
	\bproof For  \eqref{estimate1}, we fix a point $x\in M$ and assume that $h_{0,\alpha\bar\beta}(x)=\delta_{\alpha\beta}$ and $h_{\alpha\bar\beta}(x)=\lambda_\alpha \delta_{\alpha\beta}$.  In this case $H(x)=\sum \lambda_\alpha(x) e^\alpha\ts e_\alpha$ and the condition $C_4^{-1}h_0\leq h\leq C_4h_0$ gives \beq C_4^{-1}\leq \lambda_\mu\leq C_4.\eeq  If $A=A_\alpha^\beta e^\alpha\ts e_\beta$, then  $h_0(A,A)=\sum_{\alpha,\beta} |A_\alpha^\beta|^2$. Moreover, one has
	$$A\cdot H^{-1}(x)=\sum_{\alpha,\beta} A_\alpha^\beta(x)\lambda^{-1}_\beta e^\alpha\ts e_\beta$$
	and 
	\beq h_0(A\cdot H^{-1}, A)=\sum_{\alpha,\beta} |A_\alpha^\beta|^2\lambda^{-1}_\beta, \quad h(A\cdot H^{-1}, A\cdot H^{-1})=\sum_{\alpha,\beta} |A_\alpha^\beta|^2\lambda^{-1}_\alpha\lambda^{-1}_\beta. \eeq
	Therefore, 
	 \beq C^{-1}_4 h_0(A\cdot H^{-1}, A)\leq h(A\cdot H^{-1}, A\cdot H^{-1}) \leq C_4 h_0(A\cdot H^{-1}, A)\eeq 
	 and 
	  \beq C_4^{-2} h_0(A,A) \leq h(A\cdot H^{-1}, A\cdot H^{-1}) \leq C_4^2 h_0(A,A).\eeq 
	 Moreover, at the point $ x \in M $, we can also assume $ g_{i\bar j}(x) = \delta_{ij} $. Hence
	 \beq \mathrm{tr}_E\left(\Lambda_{\omega_g}\left(\sq \p^{h_0} H\cdot H^{-1}\cdot \bp H\right)\right)= \sum_{i=1}^n \mathrm{tr}_E\left(\nabla_i^{h_0}H \cdot H^{-1} \cdot \bar\p_i H\right). \eeq
	On the other hand, one has 
	\beq \mathrm{tr}_E\left(\nabla_i^{h_0}H \cdot H^{-1} \cdot \bar\p_i H\right)= h_0\left(\nabla_i^{h_0}H \cdot H^{-1}, \nabla^{h_0}_i H\right). \label{selfadjoint}\eeq 	
	Indeed, if we write $\nabla_i $ for $\nabla_i^{h_0}$, then for any $A\in \Gamma(M,E^*\ts E)$, one has
	$$  \mathrm{tr}_E(\nabla_i H \cdot A \cdot \bar\p_i H) = \left(\nabla_i H_\alpha^\beta\right) A_\beta^\gamma \left(\bar\p_i H_\gamma^\alpha \right)= h_0^{\alpha\bar p}h_{0,\gamma\bar q} \left(\nabla_i H_\alpha^\beta\right) A_\beta^\gamma \left(\overline{\nabla_i H_p^q}\right)=h_0(\nabla_i H \cdot A, \nabla_i H), $$
	where the second identity holds since 
	$ H $ is Hermitian with respect to $ h_0 $ and
	\beq \bar\p_i H_\gamma^\alpha = h_0^{\alpha\bar p}h_{0,\gamma\bar q}\left(\overline{\nabla_i H_p^q}\right). \label{Hermitianrelation} \eeq
Therefore, by \eqref{selfadjoint} and \eqref{estimate1}, one obtains
	\be \mathrm{tr}_E\left(\Lambda_{\omega_g}\left(\sq \p^{h_0} H\cdot H^{-1}\cdot \bp H\right)\right)&=&\sum_i  h_0\left(\nabla_i^{h_0}H \cdot H^{-1}, \nabla^{h_0}_i H\right)\\ &\geq &C_4^{-1}\sum_{i=1}^n|\nabla^{h_0}_i H \cdot H^{-1}|_h^2\\ &=& C_4^{-1} |\p^{h_0}H \cdot H^{-1}|_h^2.\ee
	This completes the estimate \eqref{estimate2}.
 	\eproof

	\bproposition
	\label{trEH} Let $E$ be a holomorphic vector bundle over  a compact K\"ahler manifold. Suppose that $h_0$ and $h$ are two Hermitian metrics on $E$, and $\Lambda_{\omega_g}\left(\sq R^{h_0}\right)>0$. 	
	If there exists a constant  $C_4\geq 1 $ such that 
	\beq C_4^{-1}h_0\leq h\leq C_4h_0, \eeq
	then
	\beq \Delta_{\C}\mathrm{tr}_E H  \geq  C_4^{-1}|\p^{h_0}H \cdot H^{-1}|^2_h -\mathrm{tr}_E\Phi, \eeq
	where $\Delta_\C =\mathrm{tr}_{\omega_g}\sq \p\bp$ on smooth functions.
	\eproposition
	
	\bproof A straightforward computation gives 
	\begin{eqnarray}
	\p^{*,h_0 }(\p^{h_0}H \cdot H^{-1}) \cdot H  \nonumber& = & \smo\Lambda_{\omega_g}\bar\p(\p^{h_0}H \cdot H^{-1}) \cdot H \\
	& = & \p^{*,h_0}\p^{h_0} H - \smo\Lambda_{\omega_g}(\p^{h_0}H \cdot \bar\p H^{-1}) \cdot H \label{calabic31} \\
	\nonumber & = & \Delta_{\p^{h_0}} H + \smo\Lambda_{\omega_g}(\p^{h_0}H \cdot H^{-1} \cdot \bar\p H).
	\end{eqnarray}
	On the other hand, by \eqref{Currelation}, one has
	\beq \Phi - \Omega \cdot H = \p^{*,h_0}(\p^{h_0}H \cdot H^{-1}) \cdot H. \eeq
	Moreover, $\mathrm{tr}_E H=h_0(H, \mathrm{Id}_E)$ and so
	\be \Delta_\C \mathrm{tr}_E H&=&\Lambda_{\omega_g} \sq \p\bp h_0(H, \mathrm{Id}_E)\\
	&=&-\sq \Lambda_{\omega_g} h_0(\bp \p^{h_0} H, \mathrm{Id}_E)\\
	&=&-h_0\left(\Delta_{\p^{h_0}}H, \mathrm{Id}_E\right)=-\mathrm{tr}_E\left( \Delta_{\p^{h_0}}H\right). \ee 
	By taking the trace over $E$ of the equation in \eqref{calabic31}, one obtains
	\be \Delta_{\C}\mathrm{tr}_E H  =  -\mathrm{tr}_E(\Delta_{\p^{h_0}} H) &=& \mathrm{tr}_E\left(\smo\Lambda_{\omega_g}(\p^{h_0}H \cdot H^{-1} \cdot \bar\p H)\right) - \mathrm{tr}_E(\p^{*,h_0}(\p^{h_0}H \cdot H^{-1}) \cdot H ) \\
	& \geq & C_4^{-1}|\p^{h_0}H \cdot H^{-1}|^2_h + \mathrm{tr}_E(\Omega \cdot H - \Phi) \\
	& \geq &C_4^{-1}|\p^{h_0}H \cdot H^{-1}|^2_h -\mathrm{tr}_E\Phi , \ee
	where the first inequality follows from \eqref{estimate2} and the second inequality holds since $ \Omega >0$ with respect to $h_0$.
	\eproof
	
	Let $ \sE: = T^{*1,0}M \otimes E^* \otimes E $ be the Hermitian holomorphic vector bundle with the Hermitian metric induced by $ g $ on $ T^{*1,0}M $ and $ h $ on $ E $.  Let  $ T $ be the tensor
	\beq T := \p^{h_0}H \cdot H^{-1}  \in \Omega^{1,0}(M,E^*\otimes E) \simeq \Gamma(M,\sE). \eeq
	By using similar ideas as in the higher-order estimates for the Calabi-Yau equation, we derive the following estimate: 
	\bproposition
	\label{DeltaS}Let $E$ be a holomorphic vector bundle over  a compact K\"ahler manifold. Suppose that $h_0$ and $h$ are two Hermitian metrics on $E$. 	
	If there exists a constant  $C_4\geq 1 $ such that 
	\beq C_4^{-1}h_0\leq h\leq C_4h_0, \eeq
	then
	one has 
\beq \Delta_{\C} |T|_h^2\geq - C_{10}\left(|T|^2_h+1\right) \eeq
where $C_{10}$ depends on $M,\omega_g, C_4, h_0$ and an upper bound of $|\Phi|_{C^1(M,\omega_g, h_0)}$. 
	\eproposition
	
	\bproof 
	The Bochner-Kodaira formula establishes that 
	\beq \Delta_{\C} |T|_h^2 = |\p_\sE T|_h^2 + |\bar\p T|_h^2 + 2\mathrm{Re}\left(\left\langle \Lambda_{\omega_g}\left(\smo \p_\sE\bar\p_\sE T\right), T \right\rangle_\sE \right) - \mathrm{Ric}^\sE(T,T),\label{c30} \eeq
	where $\mathrm{Ric}^\sE$ is the Hermitian-Yang-Mills tensor of $\sE$.
It is clear that 
	\beq\Lambda_{\omega_g}\left(\smo \p_\sE\bar\p_\sE T\right) = \left(g^{i\bar j} \nabla^h_i \nabla^h_{\bar j} T_{k\alpha}^\beta\right) \; dz^k \otimes e^\alpha \otimes e_\beta,\label{c31}\eeq
	where $\nabla^h$ is the Chern connection on $\sE$. Note that,  we computed in \eqref{torsion} that
	\beq T= T_{k\alpha}^\beta \; dz^k \otimes e^\alpha \otimes e_\beta=\left( (\Gamma^h)_{i\alpha}^\beta - (\Gamma^{h_0})_{i\alpha}^\beta \right) dz^i \otimes e^\alpha \otimes e_\beta.\eeq 
	Hence,
	\begin{eqnarray} -g^{i\bar j} \nabla^h_i \nabla^h_{\bar j} T_{k\alpha}^\beta \nonumber& = & g^{i\bar j}\nabla^h_i \left((R^h)_{k\bar j\alpha}{}^\beta - (R^{h_0})_{k\bar j\alpha}{}^\beta \right) \\
	& = & g^{i\bar j}\nabla^h_i (R^h)_{k\bar j\alpha}{}^\beta - g^{i\bar j}(\nabla_i^h - \nabla_i^{h_0})(R^{h_0})_{k\bar j\alpha}{}^\beta - g^{i\bar j}\nabla^{h_0}_i (R^{h_0})_{k\bar j\alpha}{}^\beta, \label{c32}
	\end{eqnarray}
	where $\nabla^{h_0}$ is the Chern connection defined by $h_0$ and $\omega_g$. 
By using the Bianchi identity, one has
\beq \nabla^{h}_i (R^{h})_{k\bar j\alpha}{}^\beta=\nabla^{h}_k (R^{h})_{i\bar j\alpha}{}^\beta, \eeq 
and so we get 
	\be -g^{i\bar j} \nabla^h_i \nabla^h_{\bar j} T_{k\alpha}^\beta 
&=&g^{i\bar j}\nabla^h_k (R^h)_{i\bar j\alpha}{}^\beta - g^{i\bar j}(\nabla_i^h - \nabla_i^{h_0})(R^{h_0})_{k\bar j\alpha}{}^\beta - g^{i\bar j}\nabla^{h_0}_k (R^{h_0})_{i\bar j\alpha}{}^\beta\\
& = & \nabla^h_k (K^h)_\alpha^\beta - g^{i\bar j} (\nabla_i^h - \nabla_i^{h_0})(R^{h_0})_{k\bar j\alpha}{}^\beta -\nabla^{h_0}_k (K^{h_0})_{\alpha}{}^\beta,
\ee
By using the identity $ K^h = \Phi \cdot H^{-1} $, we get
	\beq  \nabla^h_k (K^h)_\alpha^\beta = \nabla^h_k (\Phi \cdot H^{-1})_\alpha^\beta = \left(\nabla^h_k \Phi_\alpha^\gamma\right) (H^{-1})_\gamma^\beta + \Phi_\alpha^\gamma \nabla^h_k(H^{-1})_\gamma^\beta. \eeq 
Since $h_0$ and $h$ are equivalent,
\beq |H^{-1}|_h\leq \sqrt{r} C_4\eeq 
and 
	\beq \nabla^h_k \Phi_\alpha^\gamma = \nabla^{h_0}_k \Phi_\alpha^\gamma + (\nabla^h_k-\nabla^{h_0}_k)\Phi_\alpha^\gamma, \eeq
we conclude that 	
	\beq |\left(\nabla^h_k \Phi_\alpha^\gamma\right) (H^{-1})_\gamma^\beta|_h \leq \tilde C_5|\p^{h_0}\Phi|_h+\tilde C_6|T|_h |\Phi|_h.\eeq 
Similarly, 
	\beq \nabla_k^h(H^{-1})^\alpha_\gamma = \nabla_k^h(h^{\alpha\bar\beta}h_{0,\gamma\bar\beta}) = h^{\alpha\bar\beta} \nabla_k^h h_{0,\gamma\beta} = h^{\alpha\bar\beta} \left(\nabla_k^h - \nabla_k^{h_0}\right) h_{0,\gamma\bar \beta}. \eeq
	This implies 
	\beq |\nabla^h_k (K^h)_\alpha^\beta|_h\leq C_5|\p^{h_0}\Phi|_h+C_6|T|_h |\Phi|_h \label{c33}\eeq 
	where $\tilde C_5, \tilde C_6, $ $C_5$ and $C_6$ are constants depending on $M,\omega_g, h_0, r$ and $ C_4$. It is easy to see that 
	\beq |(\nabla_i^h - \nabla_i^{h_0})(R^{h_0})_{k\bar j\alpha}{}^\beta|_h\leq C_7 |T|_h \label{c34}\eeq
	and 
	\beq |\nabla^{h_0}_k (K^{h_0})_{\alpha}{}^\beta|_h\leq C_8 \label{c35}\eeq  
		where $C_7$ and $C_8$ are constants depending on $M,\omega_g, h_0$ and $ C_4$. 
	By \eqref{c31}, \eqref{c32}, \eqref{c33}, \eqref{c34} and \eqref{c35},  we deduce that
	\beq \left|\Lambda_{\omega_g}\left(\smo \p_\sE\bar\p_\sE T\right)\right|_h \leq C_5|\p^{h_0}\Phi|_h+C_6|T|_h |\Phi|_h+C_7 |T|_h+C_8. \label{c36}\eeq
Note that $ \sE= T^{*1,0}M \otimes E^* \otimes E $ and the Hermitian-Yang-Mills tensor $\mathrm{Ric}^\sE$ is a linear combination of the Hermitian-Yang-Mills tensors  $\mathrm{Ric(\omega_g)}$ and $K^h$, i.e.,
\beq \mathrm{Ric}^\sE=- \mathrm{Ric(\omega_g)} \ts \mathrm{Id}_{E^*\ts E}+\mathrm{Id}_{T^{*1,0}M}\ts \left(-\left(K^h\right)^t\ts \mathrm{Id}_E+\mathrm{Id}_{E^*}\ts K^h\right).\eeq 
 Since $h$ and $h_0$ are equivalent and $\Phi=K^h\cdot H$,  the following estimate holds
\beq  \mathrm{Ric}^\sE(T,T)\leq C_9 |T|^2_h \left(1 +|\Phi|_h\right),\label{c37}\eeq
where $C_9=C_9(M,\omega_g, C_4, h_0)$.
By \eqref{c30}, \eqref{c36} and \eqref{c37}, one has
 	$$\Delta_{\C} |T|_h^2 \geq -\left(2C_5|\p^{h_0}\Phi|_h|T|_h+2C_6|T|^2_h |\Phi|_h+2C_7 |T|^2_h+2C_8|T|_h+C_9|T|^2_h+C_9|T|^2_h|\Phi|_h\right).$$
 	Let $c_1$ be an upper bound of $\|\Phi\|_{C^1(M,\omega_g, h_0)}$, i.e., $\|\Phi\|_{C^1(M,\omega_g, h_0)}\leq c_1$. Since $C_4^{-1}h_0\leq h\leq C_4h_0$, 
 	$$|\p^{h_0}\Phi|_h\leq C_4 c_1,\qquad |\Phi|_h\leq  C_4c_1. $$
 	 and this implies 
 $$\Delta_{\C} |T|_h^2 \geq -\left(2C_5C_4c_1|T|_h+2C_6c_1|T|^2_h +2C_7 |T|^2_h+2C_8|T|_h+C_9|T|^2_h+C_9c_1|T|^2_h\right).$$	
 By using the  Cauchy-Schwarz inequality,  we establish	\beq \Delta_{\C} |T|_h^2\geq - C_{10}\left(|T|^2_h+1\right) \eeq
 	where $C_{10}$ depends on $M,\omega_g, C_4, h_0$ and an upper bound $c_1$ of $|\Phi|_{C^1(M,\omega_g, h_0)}$. 
	\eproof

	\bproof[Proof of Theorem \ref{ThmC1estimate}] By Proposition \ref{C0estimate}, there exists $ C_4 = C_4(M,\omega_g,h_0,C_1)  $ such that 
	\beq C_4^{-1} h_0 \leq h \leq C_4h_0.\eeq  
	Moreover, 
	\beq \mathrm{tr}_E\Phi\leq \sqrt{r} |\Phi|_{h_0} \eeq  
	By Proposition \ref{trEH} and Proposition \ref{DeltaS}, 
	\beq \Delta_{\C}\mathrm{tr}_E H  \geq  C_4^{-1}|T|^2_h -\mathrm{tr}_E\Phi\geq  C_4^{-1}|T|^2_h -\sqrt{r} |\Phi|_{h_0}, \eeq and
	\beq \Delta_{\C} |T|_h^2 \geq -C_{10}(|T|_h^2 + 1), \eeq

\noindent 	For a large $ L > 0 $ such that $LC_4^{-1} - C_{10}>1$, one has 
	\beq \Delta_{g} (|T|_h^2 + L\mathrm{tr}_E H) \geq (LC_4^{-1} - C_{10})|T|_h^2 - (LC_{11} + C_{10}), \eeq
	where $C_{11}$ is an upper bound of $\sqrt{r} |\Phi|_{h_0}$. 
At a maximal point $ p \in M $ of $ (|T|_h^2 + L\mathrm{tr}_E H) $,  we obtain the estimate
	\beq |T|_h^2(p)  \leq C_{12}, \eeq
	where $C_{12}=LC_{11}+C_{10}$
Since $\mathrm{tr}_EH \leq r C_4$, we get
	\beq (|T|_h^2 + L\mathrm{tr}_E H)(p) \leq  C_{13}. \eeq
	Therefore, for any point $ x \in M $, one has 
	\beq |T|^2_h(x) \leq (|T|_h^2 + L\mathrm{tr}_E H)(x)\leq (|T|_h^2 + L\mathrm{tr}_E H)(p) \leq C_{13}. \eeq
	Since $T=\p^{h_0} H \cdot H^{-1}$, we deduce that 
	\beq |\p^{h_0} H \cdot H^{-1}|^2_h \leq C_{13}. \eeq 
By \eqref{estimate3},  one has 
	\beq |\p^{h_0} H|_{h_0}^2 \leq C_4^2 C_{13}. \eeq  
Moreover, by the relation \eqref{Hermitianrelation}, we deduce that
\beq |\bp H|_{h_0}^2 \leq  C_{4}^2 C_{13}. \eeq  
In summary, we establish the uniform $C^1$-estimate
 $|H|_{C^1(M,\omega_g, h_0)}\leq C_{14}$ where $ C_{14}$  depends on $M,\omega_g,h_0, C_1$ and an upper bound of $ |\Phi|_{C^1(M,\omega_g, h_0)}$.
	\eproof

	\vskip 2\baselineskip
	
	\section{Proofs of main results}
	
	In this section we prove Theorem \ref{main}:
		\btheorem
		\label{SecondRicci}
	Let $ E $ be a holomorphic vector bundle over a compact K\"ahler manifold $(M,\omega_g) $. Suppose that there exists a  Hermitian metric $ h_0 $ on $E$ such that its Hermitian-Yang-Mills tensor $ \Lambda_{\omega_g}\left(\sq  R^{h_0}\right)>0 $. Then for any  positive-definite Hermitian tensor $ P\in \Gamma(M,E^*\otimes \bar E^*) $,  there exists  a unique  Hermitian metric $ h $ on $E$ such that  \beq  \Lambda_{\omega_g}\left(\sq  R^{h}\right)=P.\eeq 
	\etheorem

\bremark\label{dual} The dual equation case can be addressed via Serre duality.  Suppose that there exists a  Hermitian metric $ h_0 $ on $E$ such that $ \Lambda_{\omega_g}\left(\sq  R^{h_0}\right) <0$. Then the dual metric $h_0^*$ on $E^*$ satisfies $$\Lambda_{\omega_g}\left(\sq R^{h_0^*}\right)>0$$ as a Hermitian tensor in $\Gamma(M,E\ts \bar E)$ under the natural identification $E\cong (E^*)^*$. Hence, for any  positive-definite Hermitian tensor $P_*\in \Gamma(M,E\ts \bar E)$, there exists a unique  Hermitian metric $h_*$ on $E^*$ such that 
\beq \Lambda_{\omega_g}\left(\sq R^{h_*}\right)=P_*. \eeq 
Or equivalently,  for any  negative-definite  Hermitian tensor $P\in \Gamma(M,E^*\ts \bar E^*)$, there exists a unqiue  Hermitian metric $h$ on $E$ satisfying the dual equation:
\beq \Lambda_{\omega_g} \left(\sq R^h\right)=H\cdot P\cdot \bar {H^t} \eeq 
where $H=h\cdot h_0^{-1}\in \Gamma(M,E^*\ts E)$.
\eremark

\noindent  
We begin by proving the following result, which ensures the closeness property required in Theorem \ref{SecondRicci}.	\btheorem\label{Ricciclose} Let $(M,\omega_g)$ be a compact K\"ahler manifold and $(E,h_0)$ be a Hermitian holomorphic vector bundle with $S^{h_0}=\Lambda_{\omega_g}\left(\sq  R^{h_0}\right)>0$.  Suppose that $ \{h_m\}$ is a sequence of smooth Hermitian metrics on $E$ and   $ S^{h_m}=\Lambda_{\omega_g}\left(\sq  R^{h_m}\right) $ are the corresponding Hermitian-Yang-Mills tensors in $\Gamma(M, E^*\ts \bar E^*)$. 	If 
	\beq \lim_{m}\|S^{h_m}-P\|_{C^\infty(M,\omega_g, h_0)}=0 \eeq
	for some positive-definite  Hermitian tensor $P\in  \Gamma(M,E^* \otimes \bar E^*) $, then there exists a unique  Hermitian metric $ h $ on $E$ such that 
	\beq \Lambda_{\omega_g} \left(\sq R^h\right)=P. \eeq
	\etheorem
	
	\bproof We fix several notations in $\Gamma(M,E^*\otimes E) $:
	\beq H_m := h_m\cdot h^{-1}_{0}, \quad \Phi_m := K^{h_m}\cdot H_m=S^{h_m}\cdot h_0^{-1}, \quad  P^{h_0}: = P\cdot h_0^{-1}. \eeq Hence,  $ P^{h_0}> 0 $ with respect to $h_0$.   Since $\lim_{m}\|S^{h_m}-P\|_{C^\infty(M, \omega_g,h_0)}=0 $, one has \beq \lim_m\|\Phi_m-P^{h_0}\|_{C^\infty(M, \omega_g,h_0)}=0. \eeq In particular,  there exist uniform constants $ c_1, c_2 > 0 $ and $m_0>0$, depending on $M,\omega_g, h_0$ and $P$, such that for all $ m\geq m_0 $, the following holds
	\beq c_2^{-1}\mathrm{Id}_E \leq \Phi_m \leq c_2\mathrm{Id}_E ,\quad  \nm{\Phi_m}_{C^1(M,\omega_g, h_0)} \leq c_1. \label{close1} \eeq
	By Proposition \ref{C0estimate}, there exists $c_3=c_3(M,\omega_g, h_0, c_2)$ such that 
	\beq c^{-1}_3 \mathrm{Id}_E\leq H_m\leq c_3\mathrm{Id}_E.\label{close-1}\eeq 
	By Theorem \ref{ThmC1estimate}, there exists $ C_0 = C_0(M,\omega_g,h_0, c_1,c_2, c_3)  $ such that, 
	\beq \nm{H_m}_{C^1(M,\omega_g,h_0)} \leq C_0. \label{close2}\eeq
	By the curvature relation \eqref{Currelation}, one deduces that
	\beq \Phi_m = \Omega \cdot H_m + \p^{*,h_0}(\p^{h_0}H_m \cdot H_m^{-1}) \cdot H_m.\label{close3} \eeq
Moreover,  by using similar arguments as in the proof of  Proposition \ref{trEH}, one has
	\beq \label{ellipticHn} \Delta_{\p^{h_0}} H_m = \Phi_m - \Omega \cdot H_m - \smo\Lambda_{\omega_g}(\p^{h_0} H_m \cdot H_m^{-1} \cdot \bar\p H_m) . \eeq
 Let \beq  W_m:=\Phi_m - \Omega \cdot H_m - \smo\Lambda_{\omega_g}(\p^{h_0} H_m \cdot H_m^{-1} \cdot \bar\p H_m).\label{W_mdef}\eeq  
Hence, \beq \Delta_{\p^{h_0}} H_m=W_m.\eeq  By \eqref{close1}, \eqref{close-1} and \eqref{close2}, 
 \beq \|W_m\|_{C^0(M,\omega_g, h_0)}\leq C_1(M,\omega_g, h_0,P). \eeq 
 For a fixed large integer $ p $, there exists a constant $ C_2 = C_2(M,\omega_g,h_0,P,p) > 0 $ such that 
	\beq \nm{W_m}_{L^p(M,\omega_g, h_0)} \leq C_2. \eeq 
	By $W^{2,p}$-estimate for elliptic equations, there exists $ C_3 = C_3(M,\omega_g,h_0, p) > 0 $ such that
	\beq \nm{H_m}_{W^{2,p}(M,\omega_g, h_0)} \leq C_3\left(\nm{H_m}_{L^p(M,\omega_g,h_0)} + \nm{W_m}_{L^p(M,\omega_g,h_0)}\right). \eeq 
	In particular, $ \{ H_m \} $ is uniformly bounded in $ W^{2,p}(M, \omega_g, h_0)$.
	Since $ p $ is large enough and $M$ is compact, by the compact embedding theorem, there exist a subsequence $\{H_{m_i}\}$, $\alpha\in (0,1)$ and $H\in C^{1,\alpha}(M, \omega_g, h_0)$ such that  \beq \lim_i\| H_{m_i}-H\|_{C^{1,\alpha}(M, \omega_g, h_0)}=0.\label{Hc1alpha}\eeq Moreover,  $ H$ is $h_0$-Hermitian and positive-definite. In particular, by using \eqref{W_mdef},  one has \beq  \lim_i \|W_{m_i}- W\|_{C^{0,\alpha}(M, \omega_g, h_0)}=0\eeq  where $ W \in C^{0,\alpha}(M, \omega_g, h_0) $ is given by 
	\beq  W=\Phi - \Omega \cdot H - \smo\Lambda_{\omega_g}(\p^{h_0} H \cdot H^{-1} \cdot \bar\p H).\label{Wdef}\eeq   By Schauder's estimate, there exists $ C_4 = C_4(M,\omega_g,h_0,p,\alpha) > 0 $ such that 
	\beq \nm{H_{m_i} - H_{m_j}}_{C^{2,\alpha}(M, \omega_g, h_0)} \leq C_4\left(\nm{H_{m_i} -H_{m_j}}_{C^{0,\alpha}(M, \omega_g, h_0)} + \nm{W_{m_i} - W_{m_j}}_{C^{0,\alpha}(M, \omega_g, h_0)}\right). \eeq 
	In particular, $ \{H_{m_i}\} $ is a Cauchy sequence in $ C^{2,\alpha}(M, \omega_g, h_0)$. Moreover, by \eqref{Hc1alpha}, we deduce that $H\in C^{2,\alpha}(M, \omega_g, h_0)$ and \beq \lim_i\| H_{m_i}-H\|_{C^{2,\alpha}(M, \omega_g, h_0)}=0.\eeq Letting $ m_i \> +\infty $ in \eqref{ellipticHn}, one has the following equation 
	\beq \label{ellipticH} \Delta_{\p^{h_0}} H = W = \Phi - \Omega \cdot H - \smo\Lambda_{\omega_g}(\p^{h_0} H \cdot H^{-1} \cdot \bar\p H). \eeq
 By the standard bootstrap process, one concludes that $\{H_m\}$ converges to $H$ in $C^\infty$, i.e.,  $ H \in \Gamma(M,E^*\otimes E) $.  When $m_i\>+\infty$  in \eqref{close3},
	\beq P^{h_0}= \Omega \cdot H + \p^{*,h_0}(\p^{h_0}H \cdot H^{-1}) \cdot H. \label{close4}\eeq
	If we set $ h = H \cdot h_0 $, then by \eqref{Currelation}, we obtain the Hermitian-Yang-Mills tensor relation:
	\beq S^{h}\cdot h_0^{-1}= \Omega \cdot H + \p^*(\p^{h_0}H \cdot H^{-1}) \cdot H.\label{close5} \eeq
	By \eqref{close4} and \eqref{close5}, we  conclude that 
	\beq S^h = P. \eeq
	This completes the proof.
	\eproof

	Before proving the openness part of Theorem \ref{SecondRicci}, we fix some key notions.
	Let $E$ be a holomorphic vector bundle over a compact K\"ahler manifold $(M,\omega_g)$. The subspace of smooth Hermitian sections in $\Gamma(M,E^*\otimes \bar E^*)$ is denoted by $\mathrm{Herm}(E)$, and the subspace of smooth Hermitian metrics in $\Gamma(M,E^*\otimes \bar E^*)$ is denoted by $\mathrm{Herm}^+(E)$. There is a natural map
	$G: \mathrm{Herm}^+(E) \> \mathrm{Herm}(E) $ given by the Hermitian-Yang-Mills tensor:
	\beq G(h) = S^h= S^{h}_{\alpha\bar\beta} e^\alpha\ts \bar e^\beta \in \Gamma(M,E^*\otimes \bar E^*), \eeq
	where $S^h_{\alpha\bar\beta}=g^{i\bar j} R^h_{i\bar j \alpha\bar\beta}$.
	Fix a smooth Hermitian metric $ h_0 $ on $ E $, we define 
	\beq \mathrm{Herm}(E,h_0) := \left\{ S \in \Gamma(M,E^*\otimes E)|\ S \text{ is $h_0$-Hermitian}\right\}, \eeq 
and 
	\beq \mathrm{Herm}^+(E,h_0) := \left\{ S \in \mathrm{Herm}(E,h_0)|\ S > 0  \right\}. \eeq
	It is clear that the map $F:\mathrm{Herm}(E)\> \mathrm{Herm}(E,h_0)$ given by \beq  F(h) = h \cdot h_0^{-1} \eeq  is an isomorphism and its restriction $F:\mathrm{Herm}^+(E) \> \mathrm{Herm}^+(E,h_0)$ is also an isomorphism.
	There is an induced map $\tilde G: \mathrm{Herm}^+(E,h_0) \> \mathrm{Herm}(E,h_0)$ \quad \beq   \tilde G = F \circ G \circ F^{-1}. \eeq
	
\noindent
	Since $ \mathrm{Herm}(E,h_0) $ is a closed linear subspace of $ \Gamma(M,E^*\otimes E) $, for any $ S \in \mathrm{Herm}(E,h_0) $, there is an isomorphism 
	\beq \label{tangentiso} \mathrm{Herm}(E,h_0) \simeq T_S\mathrm{Herm}(E,h_0) \eeq 
	mapping $ \Psi \in \mathrm{Herm}(E,h_0)$ to $ [\gamma_\Psi] \in T_S\mathrm{Herm}(E,h_0) $, where $ \gamma_\Psi $ is a curve 
	$\gamma_\Psi(t) = S + t\Psi$ in $\mathrm{Herm}^+(E,h_0)$.

	\blemma\label{linearization}	
	For any $ H \in \mathrm{Herm}^+(E,h_0) $,
	\beq \tilde G(H) = \p^{*,h_0}(\p^{h_0} H \cdot H^{-1})  \cdot H + \Omega \cdot H, \eeq
	where $ \Omega =K^{h_0}\in \Gamma(M,E^*\otimes  E) $. Moreover,   the linearization map \beq  \mathscr{L}_{\mathrm{Id}_E} : \mathrm{Herm}(E,h_0) \> \mathrm{Herm}(E,h_0)\eeq  of $ \tilde G $ at the point $ \mathrm{Id}_E $
is given by
	\beq \mathscr{L}_{\mathrm{Id}_E}(\Psi) = \Delta_{\p^{h_0}} \Psi + \Omega \cdot \Psi. \eeq
	\elemma
	
	\bproof	For any $ H \in \mathrm{Herm}^+(E,h_0) $, it induces a Hermitian metric 
	\beq h = F^{-1}(H) = H \cdot h_0 \in \mathrm{Herm}^+(E). \eeq
	By \eqref{secondChernRiccidiff}, one has the difference of the Hermitian-Yang-Mills tensors:
	\beq (K^h)_\alpha^\beta - (K^{h_0})_\alpha^\beta = (\p^{*,h_0}(\p^{h_0}H \cdot H^{-1}))_\alpha^\beta. \eeq
	Hence, 
	\beq S^h_{\alpha\bar\gamma} = h_{\beta\bar\gamma}(K^{h_0})_\alpha^\beta + h_{\beta\bar\gamma}(\p^{*,h_0}(\p^{h_0}H \cdot H^{-1}))_\alpha^\beta. \eeq
	That is 
	\beq G(F^{-1}(H)) = S^h_{\alpha\bar\gamma}\; e^\alpha \otimes \bar e{}^\gamma= \left(h_{\beta\bar\gamma}(K^{h_0})_\alpha^\beta + h_{\beta\bar\gamma}(\p^{*,h_0}(\p^{h_0}H \cdot H^{-1}))_\alpha^\beta\right) \; e^\alpha \otimes \bar e{}^\gamma. \eeq
	Therefore, one has 
	\be \tilde G(H)  =  F(G(F^{-1}(H))) &=& h_0^{\delta\bar\gamma}h_{\beta\bar\gamma} \left((K^{h_0})_\alpha^\beta + (\p^*(\p^{h_0}H \cdot H^{-1}))_\alpha^\beta\right) \; e^\alpha \otimes e_\delta \\
	& = &  \p^{*,h_0}(\p^{h_0} H \cdot H^{-1})  \cdot H + \Omega \cdot H. \ee
	We compute the tangent map of $\tilde G$ at $ \mathrm{Id}_E \in \mathrm{Herm}^+(E,h_0)$. For any $ \Psi \in \mathrm{Herm}(E,h_0) $, there is a curve $ \gamma(t) = \mathrm{Id}_E + t\Psi $ in $\mathrm{Herm}^+(E,h_0)$ when $t$ is small. It is easy to see that
	\beq (\mathrm{Id}_E + t\Psi)^{-1} = \mathrm{Id}_E - t\Psi + O(t^2). \eeq
	Since $ \p^{h_0}(\mathrm{Id}_E) = 0 $, one has 
	\be \tilde G(\gamma(t)) & = & t\p^{*,h_0}((\p^{h_0} \Psi) \cdot (\mathrm{Id}_E - t\Psi)) + \Omega + t\Omega \cdot \Psi + O(t^2) \\
	& = & \Omega + t\left(\Delta_{\p^{h_0}} \Psi + \Omega \cdot \Psi \right) + O(t^2). \ee
Hence, the linearization of $ \tilde G $ at $ \mathrm{Id}_E $ is 
	\beq \mathscr{L}_{\mathrm{Id}_E}(\Psi) = \Delta_{\p^{h_0}} \Psi + \Omega \cdot \Psi. \eeq
	This completes the proof.
	\eproof
	
	\noindent The following result is useful in the $C^{k,\alpha}$ and $W^{k,p}$ estimates:
	\blemma\label{transition}  Let $S\in \mathrm{Herm}(E, h_0)$ and $\Om\in \mathrm{Herm}^+(E, h_0)$. If $\Psi\in \Gamma(M,E^*\ts E)$ satisfies 
	\beq\Delta_{\p^{h_0}} \Psi + \Omega \cdot \Psi=S, \label{LinearizationPsi}\eeq
	then $\Psi\in \mathrm{Herm}(E,h_0)$. 	
	\elemma 
	
	\bproof  For any $T\in \Gamma(M,E^*\ts E)$, we show that 
	\beq \left( \p^{h_0} T\right)^*=\bp T^*, \quad \left( \bp T\right)^*=\p^{h_0} T^*. \eeq  
	Here and in the following,  the Hermitian adjoint is taken with respect to $h_0$. Indeed,  
	\beq \p \left(h_0(Tv, w)\right)=h_0(\p^{h_0}T \cdot v +T\cdot \p^{h_0} v, w)+h_0(Tv, \bp w),\eeq for any $v, w \in \Gamma(M,E)$.
	Similarly, one has 
	\beq  \p \left(h_0(v, T^*w)\right)=h_0(\p^{h_0} v, T^*w)+h_0(v, \bp T^*\cdot w+ T^* \bp w).\eeq 
	Since $h_0(Tv, w)=h_0(v, T^*w)$, we deduce that $$h_0(\p^{h_0}T \cdot v , w)=h_0(v, \bp T^*\cdot w)$$ which implies $\left( \p^{h_0} T\right)^*=\bp T^*$. By using the expressions of  $$\p\bp(h_0(Tv,w))=\p\bp(h_0(v,T^*w)),$$ we get 
	\beq \left(\p^{h_0}\bp T \right)^*=\bp \p^{h_0} T^*, \quad \left(\bp  \p^{h_0}T\right)^*=\p^{h_0} \bp T^*.\eeq 
	In particular, 
	\beq (\Delta_\bp T)^*=\Delta_{\p^{h_0}} T^*, \quad \left(\Delta_{\p^{h_0}} T\right)^*=\Delta_{\bp} T^*. \label{Delta*}\eeq 
	Since $ \Omega $ is $ h_0 $-Hermitian,   
	\beq \label{Omega*}(\Omega\cdot \Psi)^* = \Psi^* \cdot \Omega^* = \Psi^* \cdot \Omega. \eeq
	By the Bochner-Kodaira identity on vector bundle $ E^* \otimes E $, one can see
	\beq \label{BKidentity}\Delta_{\bar\p}\Psi^* - \Delta_{\p^{h_0}}\Psi^* = \Omega \cdot \Psi^* - \Psi^* \cdot \Omega. \eeq
	Therefore, by \eqref{Delta*}, \eqref{Omega*} and \eqref{BKidentity}, we obtain that  
	\beq \label{LinearizationPsi*}\Delta_{\p^{h_0}}\Psi^* + \Omega \cdot \Psi^*=\Delta_{\bp} \Psi^*+\Psi^*\cdot \Om=(\Delta_{\p^{h_0}}\Om+\Om\cdot \Psi)^* = S^* = S. \eeq
From \eqref{LinearizationPsi} and \eqref{LinearizationPsi*}, one has 
	\beq \label{diffLinearizationPsi}\Delta_{\p^{h_0}}\left(\Psi^* - \Psi\right) + \Omega \cdot (\Psi^* - \Psi) = 0. \eeq
 Since $\Om$ is $h_0$-Hermitian and positive-definite, by pairing with $ (\Psi^* - \Psi) $ in \eqref{diffLinearizationPsi}, one deduces that $ \Psi^* - \Psi =0$ and so $\Psi$ is $h_0$-Hermitian.
	\eproof 
	
	\btheorem
	\label{SecRicciopen}
	Let $ P \in \mathrm{Herm}^+(E) $. If there exists a Hermitian metric $ h_1\in \mathrm{Herm}^+(E) $ such that 
	\beq P = G(h_1)  \in \Gamma(M,E^*\otimes \bar E^*), \eeq
	then there exist open neighborhoods $ U $ of $ P \in  \mathrm{Herm}(E) $ and $ V $ of $ h_1 \in \mathrm{Herm}^+(E) $ such that $ G: U\> V $ is an isomorphism.
	
	\etheorem
	
	\bproof We fix the center point $h_1\in \mathrm{Herm}^+(E)$ and so $F(h_1)=\mathrm{Id}_E$.
	Since $\tilde G=F\circ G\circ F^{-1}$,
\beq d\tilde G_{\mathrm{Id}_E}=(dF)_{G(h_1)}\circ (dG)_{h_1}\circ (dF^{-1})_{F(h_1)}\eeq 	
Moreover, since $ F $ is an isomorphism,  we only need to verify that  $$ \mathscr{L}_{\mathrm{Id}_E} = d\tilde G_{\mathrm{Id}_E}:\mathrm{Herm}(E,h_1) \> \mathrm{Herm}(E,h_1) $$ is an isomorphism.  By Lemma \ref{linearization}, for any $\Psi\in \mathrm{Herm}(E,h_1)$,
	\beq \mathscr{L}_{\mathrm{Id}_E}(\Psi) = \Delta_{\p^{h_1}} \Psi + \Omega_1 \cdot \Psi.\eeq	
	where $\Om_1=P\cdot h_1^{-1}\in\Gamma(M,E^*\ts E)$. Since both  $\Om_1$ and $\Psi$ are
	$h_1$-Hermitian, $\mathscr{L}_{\mathrm{Id}_E}$ is  self-adjoint with respect to $\omega_g$ and $h_1$. Since $\Om_1>0$, one can see clearly that $\mathscr{L}_{\mathrm{Id}_E}$ is injective. Moreover, since  $ M $ is compact, there exist $ m_1,m_2 > 0 $ such that \beq  m_1\cdot \mathrm{Id}_E \leq \Omega_1 \leq m_2 \cdot \mathrm{Id}_E\eeq with respect to $h_1$.  Hence, 
	\beq (\mathscr{L}_{\mathrm{Id}_E}\Psi,\Psi)_{h_1}\geq (\Omega_1\cdot \Psi,\Psi)_{h_1} \geq m_1(\Psi,\Psi)_{h_1}.\eeq 
By applying the Lax-Milgram Theorem framework together with standard elliptic regularity theory, we conclude that  $\mathscr{L}_{\mathrm{Id}_E}$ is an isomorphism. For the convenience of readers, we outline key details of the proof below. For any positive integer $m$, let 
\beq W^{m,2}(M, E^*\ts E)\eeq
be the space of $W^{m,2}$ sections of $E^*\ts E$ with respect to the metric $\omega_g$ and $h_1$.
 Consider the  pairing 
	\beq  \mathscr B: W^{1,2}(M, E^*\ts E)\times W^{1,2}(M, E^*\ts E)\> \mathbb{C} \eeq defined by
	\beq \mathscr B (\Psi,T) =(\mathscr{L}_{\mathrm{Id}_E}(\Psi), T)_{h_1}\eeq in the weak sense.  It is clear that  \beq \mathscr B(\Psi,T) = (\p^{h_1}\Psi, \p^{h_1} T) + (\Omega_1 \cdot \Psi, T)\eeq
and so
\beq |\mathscr B(\Psi,T)| \leq C_1\nm{\Psi}_{W^{1,2}} \nm{T}_{W^{1,2}}, \quad |\mathscr B(\Psi,\Psi)| \geq C_2 \nm{\Psi}_{W^{1,2}}^2, \eeq
where $C_1=C_1(M,\omega_g,h_1)>0$ and $C_2=C_2(M,\omega_g,h_1)>0$.	
  Since $\mathscr{L}_{\mathrm{Id}_E}$ is self-adjoint,	by the Lax-Milgram Theorem, one establishes an isomorphism 
	\beq \label{weakiso}\mathscr{L}_{\mathrm{Id}_E}:W^{1,2}(M, E^*\ts E)\simeq \left(W^{1,2}(M, E^*\ts E)\right) ^*. \eeq
	We claim that for all $ k \geq 1 $, 
	\beq \mathscr{L}_{\mathrm{Id}_E}: W^{k+2,2}(M, E^*\ts E) \> W^{k,2}(M, E^*\ts E)\label{Lk+2}\eeq
	are isomorphisms. We only need to show the injectivity and surjectivity.
	\bd
	\item(Injectivity) Suppose that there exists  some $ \Psi \in W^{k+2,2}(M, E^*\ts E)$ such that \beq \mathscr{L}_{\mathrm{Id}_E}(\Psi) = 0 \in W^{k,2}(M, E^*\ts E), \eeq  then by pairing with $\Psi$, one deduces $ \Psi = 0 $.
	\item(Surjectivity) There is a natural inclusion (e.g. \cite[p.~291]{Bre11})
	\beq W^{k,2}(M,E\ts E^*)\subset W^{1,2}(M,E\ts E^*) \subset  \left(W^{1,2}(M,E\ts E^*)\right)^*.\label{inclusion}\eeq
  By \eqref{weakiso},
	for each $ S \in W^{k,2}(M, E^*\ts E)$, 
there exists some \beq  \Psi \in W^{1,2}(M, E^*\ts E)\eeq  such that 
	\beq \mathscr{L}_{\mathrm{Id}_E}(\Psi)=\Delta_{\p^{h_1}} \Psi + \Omega_1 \cdot \Psi = S \eeq
	in the weak $W^{1,2}(M, E^*\ts E)$ sense. By standard regularity theory for elliptic PDEs (e.g. \cite[p.~329]{Eva10}), one concludes that $ \Psi \in W^{k+2,2}(M, E^*\ts E)$. Hence, $\mathscr{L}_{\mathrm{Id}_E}$ in \eqref{Lk+2} is surjective.
	\ed
By the Sobolev embedding theorem and Lemma \ref{transition}, we deduce that the linearization map $\mathscr{L}_{\mathrm{Id}_E} :\mathrm{Herm}(E,h_1) \> \mathrm{Herm}(E,h_1)$ is an isomorphism and so Theorem \ref{SecRicciopen} holds by the implicit function theorem.
	\eproof
	
By refining the openness and closeness arguments in the proofs of Theorem \ref{Ricciclose} and Theorem \ref{SecRicciopen}, we establish an existence result under the $C^{k,\alpha}$-norm.   For any $k>1$ and $\alpha\in (0,1)$, let 
$C^{k,\alpha}\mathrm{Herm}(E)$ and $C^{k,\alpha}\mathrm{Herm}^+(E)$ be $C^{k,\alpha}$ Hermitian sections and $C^{k,\alpha}$ Hermitian metrics on $E$ respectively. Here, the $C^{k,\alpha}$-norm is taken with respect to $\omega_g$. For any $h_1\in C^{k,\alpha}\mathrm{Herm}^+(E)$, one can also define $C^{\ell ,\alpha}\mathrm{Herm}(E,h_1)$ and $C^{\ell,\alpha}\mathrm{Herm}^+(E,h_1)$ for any non-negative integer $\ell\leq k$. Moreover, we have the following diagram, which refines the construction in the proof of Theorem \ref{SecRicciopen}:	
	$$    \CD
C^{k+2,\alpha}\mathrm{Herm}^+(E) @>G_{k,\alpha}>> C^{k,\alpha}\mathrm{Herm}(E)\\
	@V_{F_{k+2,\alpha}}  VV @V_{F_{k,\alpha}}  VV  \\
	C^{k+2,\alpha}\mathrm{Herm}^+(E,h_1) @>\tilde G_{k,\alpha}>> C^{k,\alpha}\mathrm{Herm}(E,h_1).
	\endCD
	$$
	
	\btheorem
	\label{SecondRicciCk}
	Let $ E $ be a holomorphic vector bundle over a compact K\"ahler manifold $(M,\omega_g) $. Suppose that there exists a smooth Hermitian metric $ h_0 $ on $E$ such that the Hermitian-Yang-Mills tensor $ S^{h_0} > 0 $. Then for any integer $k>1$, $\alpha\in (0,1)$ and  any  $ P \in C^{k,\alpha}{\mathrm{Herm}^+(E)} $, there exists a unique Hermitian metric $ h \in C^{k+2,\alpha}{\mathrm{Herm}^+(E)} $  such that its Hermitian-Yang-Mills tensor $S^h$  coincides with  $P$.
	\etheorem
	
	\bproof
	The uniqueness of metric $h$ is established in Corollary \ref{uniqueness} since $h$ is $C^3$. For the existence of such a metric,  we define 
	\beq \mathscr{R} = \{ P \in C^{k,\alpha}{\mathrm{Herm}^+(E)}  |\  S^h=P \text{ for some}\  h \in C^{k+2,\alpha}{\mathrm{Herm}^+(E)} \} . \eeq 
	By a similar argument as in the proof of Theorem \ref{SecRicciopen}, one deduces that  $ \mathscr{R} $ is open in $ C^{k,\alpha}{\mathrm{Herm}^+(E)}  $.  Actually, for any $h_1\in C^{k+2,\alpha}{\mathrm{Herm}^+(E)}$ with $S^{h_1}>0$,  one can show that the linearization map 
	\beq  \mathscr{L}_{\mathrm{Id}_E}: C^{k+2,\alpha}{\mathrm{Herm}(E,h_1)}  \> C^{k,\alpha}{\mathrm{Herm}(E,h_1)} \label{ckiso} \eeq 
  of $\tilde G_{k,\alpha}$ at point $\mathrm{Id}_E$ is	given by 
	\beq \mathscr{L}_{\mathrm{Id}_E}(\Psi) = \Delta_{\p^{h_1}} \Psi + \Omega_1 \cdot \Psi. \eeq 
Moreover, it	is an isomorphism. Indeed,  there is an inclusion
	\beq C^{k,\alpha}{\mathrm{Herm}(E,h_1)} \subset W^{1,2}(M,E^*\ts E) \subset  \left(W^{1,2}(M,E^*\ts E)\right)^*. \label{weakiso2}\eeq
	 By using similar argument as in the proof of Theorem \ref{SecRicciopen},
	for any $ S \in C^{k,\alpha}{\mathrm{Herm}(E,h_1)}$, 
	by \eqref{weakiso2}, there exists some \beq  \Psi \in W^{1,2}(M,E^*\ts E)\eeq  such that 
	\beq \Delta_{\p^{h_1}} \Psi + \Omega_1 \cdot \Psi=S \eeq
	in the weak $W^{1,2}(M, E^*\ts E)$ sense. By standard regularity theory for elliptic PDEs,  \beq  \Psi \in W^{k+2,2}(M,E^*\ts E).\eeq  By  the Sobolev embedding theorem, $\Psi \in W^{3,p_1}(M,E^*\ts E)$ with \beq \frac{1}{p_1}=\frac{1}{2}-\frac{k-1}{2n}\eeq  By repeating the above process, we know \beq \Psi \in W^{k+2,p_1}(M,E^*\ts E).\eeq By induction, one can show that 
	\beq \Psi \in W^{k+2,p_m}(M,E^*\ts E)\eeq
	where 
	\beq \frac{1}{p_m}=\frac{1}{p_{m-1}}-\frac{k-1}{2n}.\eeq  
	By the induction formulation,  
	\beq \frac{1}{p_m}=\frac{1}{2}-\frac{m(k-1)}{2n}.\eeq 
	Let $m_0=\lceil\frac{n}{k-1}\rceil-1$. Then 
	\beq \frac{1}{p_{m_0}}\leq \frac{k-1}{2n}\leq \frac{k-\alpha}{2n}.\eeq 
	In particular, by the Sobolev embedding theorem, we obtain
	\beq W^{k+2,p_{m_0}}(M,E^*\ts E)\subset W^{k+2,\frac{2n}{k-\alpha}}(M,E^*\ts E)\subset  C^{2,\alpha}(M,E^*\ts E).\eeq 
  Hence, $\Psi\in C^{2,\alpha}(M,E^*\ts E)$. Since $ S \in C^{k,\alpha}(M,E^*\ts E)$ and $\Om_1\in C^{k,\alpha}(M,E^*\ts E) $, 	by Schauder's theory for the elliptic equation $$\Delta_{\p^{h_1}} \Psi + \Omega_1 \cdot \Psi=S, $$ we conclude that $\Psi \in C^{k+2,\alpha}(M,E^*\ts E)$.  Hence by Lemma \ref{transition}, the map in \eqref{ckiso} is an isomorphism and so $ \mathscr{R} $ is open in $ C^{k,\alpha}{\mathrm{Herm}^+(E)}  $.
   By using similar arguments as in the proof of Theorem \ref{Ricciclose},  one can show that: 
  if  $ \{h_m\}$ is a sequence of $C^{k+2,\alpha}$-Hermitian metrics on $E$ and  
   \beq \lim_{m}\|S^{h_m}-P\|_{C^{k,\alpha}(M,\omega_g, h_0)}=0 \eeq
   for some  $P\in C^{k,\alpha}\mathrm{Herm}^+(E) $, then there exists a unique $h\in C^{k+2,\alpha}{\mathrm{Herm}^+(E)} $ such that 
   $S^h=P$.   Hence, $ \mathscr{R} $ is  closed in $ C^{k,\alpha} {\mathrm{Herm}^+(E)} $. Since $ C^{k,\alpha} {\mathrm{Herm}^+(E)} $ is connected,  one deduces that $  \mathscr{R}= C^{k,\alpha} {\mathrm{Herm}^+(E)}$.
	\eproof

	\bproof[Proof of Theorem \ref{SecondRicci}]
	 For any $ P \in \mathrm{Herm}^+(E) $, by Theorem \ref{SecondRicciCk}, there exists a unique Hermitian metric $ h_k \in C^{k+2,\alpha} {\mathrm{Herm}^+(E)} $ such that 
	\beq S^{h_k} = P. \eeq	
Moreover, by Corollary \ref{uniqueness},  $  h_k $  are the same for all $k>1$ and we denote it by $h$. Hence, $ h $ is a smooth Hermitian metric on $E$ and $ S^h = P $. 
	\eproof
	\vskip 1\baselineskip
	
\noindent 	We demonstrate that Theorem \ref{main} holds in a  more general setting. 	
Let $h$ be a smooth Hermitian metric on $ E $. The minimal eigenvalue function  $ \kappa_h: M \> \mathbb{R} $  of the Hermitian-Yang-Mills tensor $K^h$ is defined as
\beq \kappa_h(x) = \inf_{0 \neq v \in E_x} \frac{h(K^hv,v)}{h(v,v)}. \eeq
It is easy to see that $K^h > 0 $ with respect to $h$ if and only if $ \kappa_h(x) > 0 $ for any $ x \in M $.

\btheorem\label{main4}
Let $ E $ be a holomorphic vector bundle over a compact K\"ahler manifold $(M,\omega_g) $. Suppose that there exists a Hermitian metric $ h_0 $ on $E$ such that 
\beq \int_M \kappa_{h_0}(x) \omega_g^n > 0. \label{integral}\eeq Then for any Hermitian positive definite tensor $ P\in \Gamma(M,E^*\otimes \bar E^*) $,  there exists  a unique smooth Hermitian metric $ h $ on $E$ such that  \beq \Lambda_{\omega_g}\left( \sq R^h\right)=P.\eeq 
\etheorem

\bproof By using the condition \eqref{integral}, there exists  some $\lambda_0>0$ such that 
\beq \int_M \kappa_{h_0}(x) \omega_g^n = \lambda_0\int_M \omega_g^n. \eeq
Since $\kappa_{h_0}$ is Lipschitz,  by Hodge theory, there exists a function $ f \in C^{2,\alpha}(M,\mathbb{R}) $ such that 
\beq \kappa_{h_0}+ \Delta_{\C} f = \lambda_0. \eeq
We define a new metric $\tilde h = e^{-f}h_0 $.  By \eqref{secondChernRiccidiff}, one has 
\beq \Lambda_{\omega_g}\left(\sq \Theta^{\tilde h}\right) - \Lambda_{\omega_g}\left(\sq \Theta^{ h_0}\right)= \p^{*,h_0}\left( (\p^{h_0} H) \cdot H^{-1} \right) \eeq
where $H=\tilde h\cdot h_0^{-1}= e^{-f} \cdot \mathrm{Id}_E$. A straightforward computation shows 
\be \p^{*,h_0}\left( (\p^{h_0} H) \cdot H^{-1} \right) & = & \p^{*,h_0}\left( -\left(e^{-f}\p f \otimes \mathrm{Id}_E\right) \cdot H^{-1} \right) \\
& = & \smo\Lambda_{\omega_g}\bar\p \left( -\p f \otimes \mathrm{Id}_E \right) \\&=& -\smo\Lambda_{\omega_g} \bar\p\p f \cdot \mathrm{Id}_E = \Delta_{\C}f \cdot \mathrm{Id}_E. \ee
Moreover, one has 
\be \kappa_{\tilde h}(x) = \inf_{0 \neq v \in E_x} \frac{{\tilde h}(K^{\tilde h}v,v)}{{\tilde h}(v,v)} &=& \inf_{0 \neq v \in E_x} \frac{{\tilde h}(K^{h_0}v,v)}{{\tilde h}(v,v)} + \Delta_{\C}f(x) \\ 
& = & \inf_{0 \neq v \in E_x} \frac{h_0(K^{h_0}v,v)}{h_0(v,v)} + \Delta_{\C}f(x) \\&=& \kappa_{h_0}(x) + \Delta_{\C}f(x)=\lambda_0>0. \ee
In particular, the Hermitian-Yang-Mills tensor $\Lambda_{\omega_g}\left(\sq R^{\tilde h}\right)$ is positive- definite. By a perturbation trick, one can find a smooth Hermitian metric $h_1$ close to $\tilde h$ such that $\Lambda_{\omega_g}\left(\sq R^{ h_1}\right)$ is positive-definite.   Now Theorem \ref{main4} follows by  Theorem \ref{main}.
\eproof 
	
		\bremark Theorem \ref{main4} also holds on general Hermitian manifolds and  provides an alternative proof of the main results  in \cite{Yang24} and \cite{LL22}.
	\eremark
	
		\vskip 2\baselineskip

	\section{Chern number inequalities}
	
	In this section, we prove Theorem \ref{main3} and Theorem \ref{main5}.
	\btheorem Let $(M,\omega_g)$ be a compact K\"ahler manifold and $E$ be a holomorphic vector bundle of rank $r$ over $M$.  Suppose that there exists a Hermitian metric $h$ on $E$ satisfying 
	\beq  a\cdot h \leq  \Lambda_{\omega_g} \left(\smo  R^{h}\right) \leq b \cdot h, \eeq 
	for some constants $a, b\in \R$. Then the following Chern number inequality holds
\beq \int_M \left( (r-1)c^2_1(E)-2rc_2(E) \right) \wedge \omega_g^{n-2} \leq \frac{r(r-1)\left(b - a\right)^2}{8\pi^2n^2}\int_M \omega_g^n.\label{CN} \eeq
\etheorem
	\bproof  
  	In the following, we work on the abstract Hermitian vector bundle $(E,h)$.
	  Let $\{z^i\}$ be local holomorphic coordinates on $M$ and $\{e_\alpha\}$ be holomorphic frames of $E$.  The Chern curvature  tensor $R^h\in \Gamma(M,\Lambda^{1,1}T^*M\ts E^*\ts \bar E^*)$ is given by 
	$$R^h=R^h_{i\bar j \alpha\bar \beta} dz^i\wedge d\bar z^j\ts  e^\alpha\ts \bar e^\beta.$$
	We set $$R^{(1)}_{i\bar j} =h^{\alpha\bar\beta}R^h_{i\bar j \alpha\bar \beta}, \quad R^{(2)}_{\alpha\bar\beta}=g^{i\bar j}  R^h_{i\bar j \alpha\bar\beta } \qtq{and} s_h =g^{i\bar j} R^{(1)}_{i\bar j}.$$
   The following expressions are well-known:
	\beq \label{c1computation} \int_M c^2_1(E) \wedge \omega_g^{n-2} = \frac{1}{4\pi^2n(n-1)}\int_M \left(s_h^2 - \left|\mathrm{Ric}^{(1)}\right|_g^2 \right) \omega_g^n, \eeq
	and 
	\beq \label{c2computation} \int_M c_2(E) \wedge \omega_g^{n-2} = \frac{1}{8\pi^2n(n-1)}\int_M \left(s_h^2 - \left|\mathrm{Ric}^{(1)}\right|_g^2 - \left| \mathrm{Ric}^{(2)}\right|^2_h + |R^h|_{g,h}^2\right) \omega_g^n. \eeq
	We define $T\in \Gamma(M,\Lambda^{1,1}T^*M\ts E^*\ts \bar E^*)$:
	\beq T_{i\bar j \alpha\bar\beta } = R^h_{i\bar j\alpha\bar\beta} - \frac{1}{n}g_{i\bar j}R^{(2)}_{\alpha\bar\beta} - \frac{1}{r}R^{(1)}_{i\bar j}h_{\alpha \bar\beta} + \frac{1}{nr}g_{i\bar j}h_{\alpha \bar\beta}s_h. \eeq
	A straightforward computation shows that 
	\beq |T|^2=|R^h|^2 - \frac{1}{n}\left|\mathrm{Ric}^{(2)}\right|^2 - \frac{1}{r}\left|\mathrm{Ric}^{(1)}\right|^2 + \frac{1}{nr} s_h^2.\eeq 
	By using $|T|^2$, one has
	\be && \int_M \left(2rc_2(E) - (r-1)c^2_1(E)\right) \wedge \omega_g^{n-2} \\
	& = & \frac{1}{4\pi^2n(n-1)}\int_M \left(s_h^2 - \left|\mathrm{Ric}^{(1)}\right|^2 - r\left|\mathrm{Ric}^{(2)}\right|^2 + r|R^h|^2\right) \omega_g^n \\
	& = & \frac{1}{4\pi^2n(n-1)}\int_M \left(r|T|^2 + \frac{n-1}{n}\left( s_h^2 - r\left|\mathrm{Ric}^{(2)}\right|^2\right) \right) \omega_g^n. \ee
		Let $\{\lambda_1,\cdots, \lambda_r\}$ be the eigenvalues of $ \Lambda_{\omega_g}\sq R^h $  with respect to $ h $.  Since $$ a h\leq \Lambda_{\omega_g}\left(\sq R^h\right) \leq bh, $$ it is obvious that 
	\beq a \leq \lambda_k \leq b, \eeq
	and so 
	\beq s_h^2 - r\left|\mathrm{Ric}^{(2)}\right|^2 = \left(\sum_{k=1}^r \lambda_k\right)^2 - r\sum_{k=1}^r \lambda_k^2 =  -\sum_{i<j} (\lambda_i-\lambda_j)^2 \geq -\frac{r(r-1)}{2}(a-b)^2. \eeq
Finally,  we establish the Chern number inequality
$$ \int_M \left(2rc_2(E) - (r-1)c^2_1(E)\right) \wedge \omega_g^{n-2}\geq  -\frac{r(r-1)\left(b - a\right)^2}{8\pi^2n^2}\int_M \omega_g^n+\frac{r}{4\pi^2n(n-1)}\int_M |T|^2\omega_g^n.$$
Hence, we obtain the result.
\eproof 
	
	The following result refines and extends Yau's classical Chern number inequalities \cite{Yau77} for  K\"ahler-Einstein manifolds.
	\btheorem
	\label{Chernnumberindentity}
	Let $ (M,\omega_g) $ be a compact K\"ahler manifold with $\dim M\geq 2$.  If there exist two constants $ a,b \in \mathbb{R} $ such that
	\beq a\omega_g \leq \mathrm{Ric}(\omega_g) \leq b\omega_g, \eeq
	then the following Chern number inquality holds 
	\beq \int_M \left(nc^2_1(M)-2(n+1)c_2(M) \right) \wedge \omega_g^{n-2} \leq \frac{(a - b)^2(n^2-2)}{8\pi^2n^2}\int_M \omega_g^n. \label{CN2}\eeq
	\etheorem
	
	\bproof
	It is known that 
	\beq \label{c11computation} \int_M c^2_1(M) \wedge \omega_g^{n-2} = \frac{1}{4\pi^2n(n-1)}\int_M \left(s^2 - \left|\mathrm{Ric}\right|^2 \right) \omega_g^n, \eeq
	and 
	\beq \label{c22computation} \int_M c_2(M) \wedge \omega_g^{n-2} = \frac{1}{8\pi^2n(n-1)}\int_M \left(s^2 - 2\left|\mathrm{Ric}\right|^2 + |R|^2\right) \omega_g^n, \eeq
	where $R$ is the curvature tensor, $\mathrm{Ric}$ is the Ricci curvature, and  $s$ is the scalar curvature of the K\"ahler metric $g$.
	Let   $ T \in \Omega^{1,1}(M,T^{*1,0}M \otimes T^{*0,1}M) $  be the tensor given by
	\beq T_{i\bar j k\bar\ell} = R_{i\bar jk\bar\ell} - \frac{1}{n}g_{i\bar j}R_{k\bar\ell} - \frac{1}{n}R_{i\bar j}g_{k\bar\ell} - \frac{1}{n(n+1)} g_{i\bar\ell}g_{k\bar j}s+ \frac{n+2}{n^2(n+1)}g_{i\bar j}g_{k\bar\ell} s . \eeq
	A straightforward computation shows that
	\beq \label{normT} |T|^2 = |R|^2 - \frac{2}{n}\left|\mathrm{Ric}\right|^2 + \frac{2}{n^2(n+1)} s^2. \eeq
	By \eqref{c11computation}, \eqref{c22computation} and \eqref{normT}, one has 
	\be && \int_M \left(2(n+1)c_2(M) - nc_1(M)^2\right) \wedge \omega_g^{n-2} \\
	& = & \frac{1}{4\pi^2n(n-1)}\int_M \left(s^2 - (n+2)\left|\mathrm{Ric}\right|^2 + (n+1)|R|^2\right) \omega_g^n \\
	& = & \frac{1}{4\pi^2n(n-1)}\int_M \left((n+1)|T|^2 + \left(1-\frac{2}{n^2}\right)\left(s^2 - n\left|\mathrm{Ric}\right|^2\right) \right) \omega_g^n. \ee
	Let $\{\lambda_1,\cdots, \lambda_n\}$ be the eigenvalues of $ \mathrm{Ric}(\omega_g) $  with respect to $ g $.  Since $ a\omega_g \leq \mathrm{Ric}(\omega_g) \leq b\omega_g $, it is obvious that $ a \leq \lambda_k \leq b$,
	and so 
	\beq s^2 - n\left|\mathrm{Ric}\right|^2 = \left(\sum_{k=1}^n \lambda_k\right)^2 - n\sum_{k=1}^n \lambda_k^2 =  -\sum_{i<j} (\lambda_i-\lambda_j)^2 \geq -\frac{n(n-1)}{2}(a-b)^2. \eeq
	Therefore, we obtain 
	\be \int_M \left(2(n+1)c_2(M) - nc_1(M)^2\right) \wedge \omega_g^{n-2} & \geq & \frac{n^2-2}{4\pi^2n^3(n-1)}\int_M \left(s^2 - n\left|\mathrm{Ric}\right|^2 \right) \omega_g^n \\
	& \geq & -\frac{(a - b)^2(n^2-2)}{8\pi^2n^2}\int_M \omega_g^n.
	\ee
	This completes the proof.
	\eproof

\bremark When $h$ is another K\"ahler or balanced metric on $M$-- i.e.,  $d\omega^{n-1}_h=0$-- and satisfies the inequality 
\beq  a\cdot h \leq  \Lambda_{\omega_g} \left( \smo  R^{h} \right)\leq b \cdot h, \eeq 
for some constants $a, b\in \R$,  then analogous Chern number inequalities to those in \eqref{CN2} hold.
\eremark

\vskip 2\baselineskip
	
	\section{The prescribed Hermitian-Yang-Mills tensor equation on line bundles}
	
	In this section, we investigate the prescribed Hermitian-Yang-Mills equation on line bundles and clarify its connection with  RC-positivity framework. Moreover, we  show that the condition $P>0$ in Theorem \ref{main} is necessary:
	\bproposition \label{counterexample}
	Let $ L $ be a holomorphic line bundle over a compact K\"ahler manifold $(M,\omega_g) $. Suppose that there exists a Hermitian metric $ h_0 $ on $ L $ such that $\Lambda_{\omega_g}\left(\sq R^{h_0}\right)>0$. If $P\in \Gamma(M,L^*\ts \bar L^*)$ fails to be Hermitian positive definite, then the equation
	\beq \Lambda_{\omega_g}\left(\sq R^{h}\right)=P \eeq 
	may admit two distinct smooth solutions.
	\eproposition
	
	\blemma Let $L$ be a holomorphic line bundle over a compact K\"ahler manifold $(M,\omega_g)$. Then for any Hermitian metric $h_0$ on $L$, 
	\beq S^{h_0}=\Lambda_{\omega_g}\left(\sq R^{h_0}\right)= s_{h_0}\cdot  h_0 \in \Gamma(M,L^*\ts \bar L^*)\eeq 
	where $s_{h_0}$ is the Chern scalar curvature of $(L,h_0)$.
	\elemma 
	\bproof 
	
	It is well-known for any Hermitian metrics $h_1, h_2 \in \Gamma(M,L^*\ts \bar L^*)$, there exists some smooth function $\phi\in C^\infty(M,\R)$ such that 
	\beq h_2=e^{-\phi}h_1.\eeq 
	Moreover,  if  $P\in \Gamma(M,L^*\ts \bar L^*)$ is a Hermitian tensor, then $P= f h_1$ some $f\in C^\infty(M,\R)$.	Let $e$ be a local holomorphic base of $L$ and $h_0=h_0(e,e)$.	 By using the curvature formula, one has 
	\beq S^{h_0}=\Lambda_{\omega_g}\sq R^{h_0}=-g^{i\bar j}\frac{\p^2\log h_0}{\p z^i\p\bar z^j} \cdot h_0 e^*\ts \bar e^*=-\left(\Delta_g\log h_0\right) h_0 e^*\ts \bar e^*.\eeq 
	The Chern scalar curvature function is 
	\beq s_{h_0}= \mathrm{tr}_{h_0}S^{h_0}=-\Delta_g\log h_0.\eeq 
	Hence, $S^{h_0}=\Lambda_{\omega_g}\left(\sq R^{h_0}\right)= s_{h_0}\cdot  h_0 $.
	\eproof

	\bproposition\label{linebundle} Let $(L,h_0)$ be a Hermitian holomorphic line bundle over a compact K\"ahler manifold $(M,\omega_g)$ and $S^{h_0}=F h_0$ for some $F\in C^\infty(M,\R)$.  Suppose that $P= Gh_0\in \Gamma(M,L^*\ts \bar L^*)$ for some $G\in C^\infty(M,\R)$. Then there exists  some Hermitian metric $h$ on $L$ satisfying
	\beq \Lambda_{\omega_g}\left(\sq R^{h}\right)=P, \eeq 
	if and only if  there exists some $\phi\in C^\infty(M,\R)$ such that 
	\beq \left(\Delta_g \phi+F\right) e^{-\phi}=G.\eeq  In this correspondence,  one has  $h=e^{-\phi} h_0$.
	\eproposition 
	
	\bproof 
	It is known that $S^{h_0}= F h_0$ where
	\beq F=-\Delta_g\log h_0.\eeq 
	Hence, if $S^{h}=P=Gh_0$ for some Hermitian metric $h$, then 
	\beq G=-\Delta_g \log h\cdot \left(h\cdot h_0^{-1}\right).\eeq  	
	On the other hand,  $h=e^{-\phi} h_0$ for some $\phi\in C^\infty(M,\R)$.  Therefore, $$-\Delta_g \log h=-\Delta_g\log h_0+\Delta_g\phi$$ and so $G=\left(\Delta_g \phi+F\right) e^{-\phi}$.
	\eproof

	\noindent 
	The following result is a consequence of Theorem \ref{main}, Theorem \ref{main4} and Proposition \ref{linebundle}, which is a line bundle version of the Kazdan-Warner theorem (\cite{KW74}).
	\bcorollary
	Let $ L $ be a holomorphic line bundle over a compact K\"ahler manifold $(M,\omega_g) $. Suppose that there exists a Hermitian metric $ h_0 $ on $ L $ such that its Chern scalar curvature $ F = s_{h_0} $ satisfying $$\int_M F\omega_g^n>0,$$ then for any positive smooth function $ G $, the following equation has the unique smooth solution $ \phi $:
	\beq e^{-\phi}(F+\Delta_{g}\phi) = G. \label{PDEsolution}\eeq
	\ecorollary
	 For a detailed discussion on the existence of Hermitian metrics with positive Chern scalar curvature,  positive total Chern scalar curvature and RC-positivity for line bundles, we refer to  \cite{Yang18}, \cite{Yang19a} and \cite{Yang19b}. This is also a motivation of this project for vector bundles. 	 In this following, we show that the equation \eqref{PDEsolution} may admit multiple solutions when $G$ is not positive.
	
	\bproof[Proof of Proposition \ref{counterexample}] 	 Let $ L $ be an ample line bundle or a line bundle with 
	$$\int_M c_1(L)\wedge \omega_g^{n-1}>0$$
	 over a projective K\"ahler manifold $ (M,\omega_g) $.  There exists a Hermitian metric $h_0$ on $L$ such that  $S^{h_0}=\Lambda_{\omega_g}\left(\sq R^h\right)=Fh_0>0$ (\cite{Yang19a}). Recall that \beq F=-\Delta_g \log h_0=s_{h_0}.\eeq 
	Let $ f \in C^\infty(M,\mathbb{R}) $ be a smooth function satisfying  
	\bd
	\item $ f\geq 0 $ and there exists $ p \in M $ such that $ f(p) = 0 $ and $ f(x) > 0 $ if $ x \neq p $.
	\item There exists some small $ \delta > 0 $ and $ q \in B(p,\delta) $ such that
	\beq f(q) = \begin{cases} e^{-1/r^n} & \text{if } q \neq p \\
		0 & \text{if } q = p,
	\end{cases} \eeq 
	where $ r = d_g(p,q) $ is the distance function. 
	\ed
	We consider the map $ Q:(0,+\infty) \> \mathbb{R} $ defined by 
	\beq Q(t) = \int_M \left(\Delta_{g}f\right)/(e^{f+t} - 1) \cdot \omega_g^n. \eeq
	We claim that 
	\beq \lim_{t\>0^+} Q(t) = +\infty, \quad \lim_{t\> +\infty} Q(t) = 0. \eeq
	Indeed,	by Stokes' theorem, one has
	\be 2Q(t) & = & 2\int_M \left(\Delta_{g}f\right)/(e^{f+t} - 1) \cdot \omega_g^n = -\int_M g\left( \nabla f,  \nabla (e^{f+t} - 1)^{-1}\right) \cdot \omega_g^n \\
	& = & \int_M |\nabla f|^2 \cdot (e^{f+t}-1)^{-2} \cdot e^{f+t} \cdot \omega_g^n = \int_M |\nabla f|^2 \cdot \left(e^{(f+t)}+e^{-(f+t)}-2\right)^{-1} \cdot \omega_g^n. \ee
	Therefore, $ Q $ is a decreasing in $t$ and by monotone convergence theorem, one has 
	\beq \quad \lim_{t\> +\infty} Q(t) = 0. \eeq
	On the other hand, by monotone convergence theorem
	\beq \lim_{t\>0^+} Q(t) = \int_{M\setminus\{p\}} |\nabla f|^2 \cdot (e^f + e^{-f} -2)^{-1} \cdot \omega_g^n. \eeq
	Moreover,	there is a constant $ A = A(M,\delta) > 0 $ such that for any $ q \in B(p,\delta) $, one has
	\beq e^{f}+e^{-f}-2 \leq Af^2, \eeq
	and so 
	\be \lim_{t\>0^+} Q(t) & \geq & A^{-1}\int_{M\setminus\{p\}} |\nabla f|^2/f^2 \cdot \omega_g^n = A^{-1}\int_{B(p,\delta) \setminus \{p\}} |\nabla \mathrm{log}(f)|^2 \cdot \omega_g^n \\
	& = & n^2A^{-1}\int_{B(p,\delta) \setminus \{p\}} r^{-(2n+2)} \cdot \omega_g^n = +\infty. \ee

 On the other hand, by continuity of $Q(t)$, there exists some $ t_0 \in (0,+\infty) $ such that
	\beq  \int_M \left(\Delta_{g} f\right)/\left(e^{f+t_0} - 1\right) \cdot \omega_g^n =2\pi n \int_M c_1(L) \wedge \omega_g^{n-1}>0. \eeq
	Let  $ \Phi = f+t_0 $. Then one has 
	\beq \label{Intphi} \int_M \left(\Delta_{g}  \Phi\right)/\left(e^{ \Phi} - 1\right) \cdot \omega_g^n=2\pi n \int_M c_1(L) \wedge \omega_g^{n-1} = \int_M s_{h_0} \cdot \omega_g^n=\int_M F \cdot \omega_g^n. \eeq
	By Hodge theory, there exists a smooth function $ \psi \in C^\infty(M,\mathbb{R}) $ such that 
	\beq \frac{\Delta_{g} \Phi}{e^{ \Phi} - 1} = F - \Delta_{g}\psi. \eeq
	If we define
	\beq G:=\frac{\Delta_g  \Phi}{e^{ \Phi}-1}\cdot e^\psi \eeq then it is obvious that
	the equation 
	\beq e^{-\phi}(F+\Delta_{g}\phi) = G \eeq
	has two different solutions $\phi=-\psi$ and $\phi=\Phi-\psi$.  In particular,  the equation 
	\beq \Lambda_{\omega_g}\sq R^{h}=Gh_0\eeq 
	has	two different solutions  $h_1= e^\psi \cdot h_0$ and  $h_2=e^{\psi- \Phi} h_0$.  	\eproof

	\bexample   Let $ (L,h_0)$ be a  Hermitian holomorphic trivial line bundle over a compact K\"ahler manifold $(M,\omega_g)$.
If  $P=G h_0\in \Gamma(M,L^*\ts \bar L^*)$ for some $G\in C^\infty(M,\R)$ satisfies 
	$$\int_M G \omega_g^n< 0,$$
	then  the equation  $\Lambda_{\omega_g}\left(\sq R^{h}\right)=P$ has no solution.
	\eexample

\end{document}